\documentclass[10pt]{amsart}

\usepackage{amssymb} 
\usepackage[mathscr]{eucal} 
\usepackage[matrix,arrow,tips,curve,ps]{xy}
\usepackage{amsmath} 
\usepackage{epsfig}
\usepackage{amscd}
\newtheorem{TEO}{Theorem}[section]
\newtheorem{PROP}[TEO]{Proposition}
\newtheorem{LEM}[TEO]{Lemma}

\newtheorem{COR}[TEO]{Corollary}
\newtheorem{REM}[TEO]{Remark}

\theoremstyle{definition}

\def\OO{{\mathcal O}}

\newcommand\de{\delta}
\newcommand\dual{\mathrel{\raise3pt\hbox{$\underline{\mathrm{\thinspace d
\thinspace}}$}}}

\newcommand\proj{\mathbb P}

\newcommand\R{\mathbb R} 
\newcommand\Co{\mathbb C}

\newlength{\spazio}

\def\R{{\mathbb R}}

\def\Co{{\mathbb C}}
\begin{document}

\footnote{
Partially supported by  FAR 2010 (PV) ``Variet\`a algebriche, calcolo algebrico, grafi orientati e topologici'', by INdAM (GNSAGA), and by MIUR of the Italian Government in the framework of the National Research Project ``Geometria algebrica e aritmetica, teorie coomologiche e teoria dei motivi''(PRIN 2008).
 }

\title{Prym map and second gaussian map for Prym-canonical line bundles}

\author[E. Colombo]{Elisabetta Colombo}
\address{Dipartimento di Matematica,
Universit\`a di Milano, via Saldini 50,
     I-20133, Milano, Italy } \email{{\tt
elisabetta.colombo@unimi.it}}

\author[P. Frediani]{Paola Frediani}
\address{ Dipartimento di Matematica, Universit\`a di Pavia,
via Ferrata 1, I-27100 Pavia, Italy } \email{{\tt
paola.frediani@unipv.it}}

\maketitle

\setlength{\parskip}{.1 in}

\begin{abstract}
We show that the second fundamental form of the Prym map lifts the
second gaussian map $\mu_A$ of the Prym-canonical bundle. We
prove, by degeneration to binary curves,  that $\mu_A$ is
surjective for the general point $[C,A]$ of ${\mathcal R}_g$ for
$g\geq 20$.
\end{abstract}

\section{Introduction}

Similarly to the period map $P_g:{\mathcal M}_g \rightarrow
{\mathcal A}_g$, the Prym map $Pr_g: {\mathcal R}_g \rightarrow
{\mathcal A}_{g-1}$ provides a way to link the geometry of moduli
spaces of curves  to the geometry of  moduli spaces of principally
polarized abelian varieties. Recall that ${\mathcal R}_g$ denotes
the moduli space which parametrizes isomorphism classes of pairs
$[C, A]$, where $C$ is a smooth curve of genus $g$ and $A \in
Pic^0(C)[2] - \{\mathcal O_C\}$ is a torsion point of order 2, or equivalently
isomorphism classes of unramified double coverings $\pi: \tilde{C}
\rightarrow C$. The Prym map associates to a point $[(C,A)] \in
{\mathcal R}_g$ the isomorphism class of the connected component of zero,
$P(C,A)$ of the kernel of the norm map $Nm_{\pi}:J{\tilde
C} \rightarrow JC$, with its principal polarization. Both the
period map and the Prym map have been extensively studied since a
long time, but also recently there have been important
developements on the birational geometry of ${\mathcal R}_g$
(\cite{fl}).

In this paper we focus on the study of the second fundamental form
of the Prym map analogously to what it has been done for the
second fundamental form of the period map and its link with the
second gaussian map. In fact in \cite{cpt} it is shown that the
second fundamental form of the period map lifts the second
gaussian map of the canonical line bundle, as stated in an
unpublished paper by Green and Griffiths (cf. \cite{green}). With
this geometrical motivation, in \cite{cf2} we investigated
curvature properties of ${\mathcal M}_g$ endowed with the Siegel
metric. In fact, we computed the holomorphic sectional curvature
of  ${\mathcal M}_g$ along the tangent directions given by the
Schiffer variations in terms of the second gaussian map. This also
suggested that the second gaussian map itself could give
interesting information on the geometry of the  curves, hence its
rank properties have been investigated in a series of papers (see
\cite{ccm}, \cite{cf1}, \cite{cf3}, \cite{cfp}).

Here we first generalize the lifting result of \cite{cpt} to the Prym map, namely at 
 the point $[C,A] \in{\mathcal R}_g $ we have the following commutative diagram

 $$
 \xymatrix{
 & {\mathcal N^{*}_{{\mathcal R}^0_g/{\mathcal A}_{g-1},[C,A]}} \cong I_2(K_C \otimes A)\ar[d]^{-\frac{1}{2 \pi i} \mu_A} \ar[r]^{II     } & S^2 \Omega^{1}_{{\mathcal R}_g^0, [C,A]} \cong S^2H^0(K_C^{\otimes 2}) \ar[dl]^{\bf{m}}\\
&H^0(K_C^{\otimes 4})
 & }
$$
where $I_2(K_C \otimes A)$ is the kernel of
the multiplication map $S^2H^0(K_C \otimes A) \rightarrow
H^0(K_C^{\otimes 2})$, $II$ is the second fundamental form of the Prym map, $m$ is the multiplication map and $\mu_A$ is the second gaussian map associated to the Prym canonical
bundle $K_C \otimes A$. This also allows us to generalize the
results of \cite{cf2} on the holomorphic sectional curvature of
${\mathcal R}_g$ with the Siegel metric induced by ${\mathcal
A}_{g-1}$ via the Prym map.

In the second part of the paper we concentrate on the study of the
second gaussian map $\mu_A$, $A \in Pic^0(C)[2]$ non trivial. The main result
is the proof of the surjectivity of $\mu_A$ for the general curve
$[C,A] \in {\mathcal R}_g$ of genus $g \geq 20$, generalizing
analogous results on the surjectivity of the second gaussian map
of the canonical line bundle for the general curve in ${\mathcal
M}_g$ for $g \geq 18$. For the canonical line bundle this
surjectivity for general curves of high genus was proved in
\cite{cf3} using curves on K3 surfaces, then the sharp result for
genus $\geq 18$ has been shown in \cite{ccm} using degeneration to
binary curves, i.e. stable curves which are the union of two
rational curves meeting transversally at $g+1$ points. Here we
generalize these degeneration techniques
  to prove the surjectivity of second gaussian maps
$\mu_A$, for the Prym-canonical bundles $K_C \otimes A$.

In particular, this shows that the locus of curves $[C,A] \in
{\mathcal R}_g$ ($g \geq 20$) for which the map $\mu_A$ is not
surjective is a proper subscheme of ${\mathcal R}_g$ and one
observes that for $g =20$ it is an effective  divisor in
${\mathcal R}_{20}$ of which we compute the cohomology class both
in ${\mathcal R}_{20}$ and in a partial compactification
$\tilde{{\mathcal R}}_{20}$ following computations developed in
\cite{fl}.

The paper is organized as follows: in Section 2 we describe the
second fundamental form, we prove that it is a lifting of the
second gaussian map $\mu_A$ and we compute the holomorphic
sectional curvature along the Shiffer variations. In Section 3 we
construct the Prym-canonical binary curves that we use for the
degeneration. In Section 4 we explicitly describe the ideal of the
quadrics containing the Prym-canonical binary curve. In Section 5
we prove, by induction on the genus, the surjectivity of $\mu_A$.
In Section 6 we compute the cohomology class  in $\tilde{{\mathcal
R}}_{20}$ of the degeneracy locus of the second gaussian map.
Finally in the Appendix we list the Maple scripts used in the
computations.

Finally we observe that the results and the techniques of section 2 are of different nature from the rest of the paper, which can be
read separately, once one has looked at the definition of the gaussian maps in subsection 2.2.

We also notice that in the proof of the surjectivity of $\mu_A$, in particular in sections 4 and 5, we follow the lines of the proof of the surjectivity of the 2nd gaussian map of the canonical line bundle given in \cite{ccm}.

\section{The 2nd fundamental form and the 2nd gaussian map}
\subsection{The second fundamental form of the Prym map}
We start by recalling the definition of  the Prym map
$$Pr: {\mathcal R}_g \rightarrow {\mathcal A}_{g-1},$$
which associates to a point $[(C,A)] \in {\mathcal R}_g$ its Prym variety $P(C,A)$ with its principal polarization.
If $\pi: \tilde{C} \rightarrow C$ is the unramified double covering associated to the pair $(C,A)$,
the Prym variety  of the double covering is the principally polarized abelian variety of dimension $g-1$
defined as the connected component of zero of the kernel of the norm map $Nm_{\pi}:J{\tilde C} \rightarrow JC$,
$$P(C,A) = Ker(Nm_{\pi})^0 \subset J{\tilde C}.$$
We recall that the Prym map is generically an embedding for $g \geq 7$ (\cite{fs}, \cite{ka}). Hence there exists
an open set ${\mathcal R}^0_g \subset {\mathcal R_g}$ where $Pr$ is an embedding and such that there exists the
universal family $f: {\mathcal X} \rightarrow {\mathcal R}^0_g$.  If $b \in {\mathcal R}^0_g$, we have
$f^{-1}(b) = (C_b, A_b)$ where $C_b$ is a smooth irreducible curve of genus $g$ and $A_b \in Pic^0(C_b)[2]$ is a
line bundle of order $2$ on $C_b$. Denote by ${\mathcal P} \in Pic({\mathcal X})$ the corresponding Prym
bundle and by ${\mathcal F}^{Pr}:= f_*(\omega_{{\mathcal X}/{\mathcal R}^0_g} \otimes {\mathcal P})$. Observe
that ${\mathcal F}^{Pr}$ is the pullback of the Hodge bundle on ${\mathcal A}_{g-1}$ to ${\mathcal R}^0_g$. More precisely, if $\psi:  {\mathcal Pr}({\mathcal X}) \to {\mathcal R}^0_g$ is the universal family of Prym varieties, so $\psi^{-1}(b)=P(C_b,A_b)$ is the Prym variety associated to the pair $(C_b, A_b)$, then $ {\mathcal F}^{Pr}$ is the Hodge bundle ${\mathcal H}^{1,0} \subset R^1 \psi_*{\Co}$ of the family $\psi$.

 On the local system $R^1 \psi_*{\Co}$ we have the flat Gauss-Manin connection $\nabla^{GM}$, and a non
 degenerate bilinear form defined as follows. At the point $P(C,A)$, the fiber of $R^1 \psi_*{\Co}$ is isomorphic to
 the vector space $H^1(\tilde{C}, {\Co})^{-}$, where $\tilde{C} \stackrel{\pi} \rightarrow C$ is the double covering
 associated to $(C,A)$ and $H^1(\tilde{C}, {\Co})^{-}$ is the anti-invariant part of the cohomology under the covering
 involution on $\tilde{C}$. If $[\omega_1], [\omega_2] \in H^1(\tilde{C}, {\Co})^{-}$, we have the following non
 degenerate bilinear form  $\langle [\omega_1],[\omega_2] \rangle = i \int_{\tilde C} \omega_1  \wedge \overline{\omega_2}$
  and the Gauss-Manin connection $\nabla^{GM}$ is compatible with it, so it induces a metric connection
  ${\nabla}^{1,0}$ on ${\mathcal H}^{1,0}$, hence a connection on ${\mathcal F}^{Pr}$ (still denoted by
  ${\nabla}^{1,0}$) and on its second symmetric power, $S^2 {\mathcal
  F}^{Pr}$. Observe that this metric on $S^2 {\mathcal
  F}^{Pr}$ is the pullback via the Prym map of the metric on
  ${\mathcal A}_{g-1}$ induced by the unique (up to scalar) $Sp(2g-2,\R)$-invariant metric on the Siegel space $H_{g-1}$.
  So we will call this metric the Siegel metric.

Consider the tangent bundle exact sequence of the Prym map
\begin{equation}
\label{tangent}
0 \rightarrow T_{{\mathcal R}^0_g} \rightarrow T_{{{\mathcal A}_{g-1}}_{|{\mathcal R}^0_g}} \rightarrow {\mathcal N_{{\mathcal R}^0_g/{\mathcal A}_{g-1}}} \rightarrow 0
\end{equation}
Its dual becomes
\begin{equation}
\label{I2}
0 \rightarrow {\mathcal I}_2 \stackrel{i} \rightarrow S^2f_*(\omega_{{\mathcal X}/{\mathcal R}^0_g} \otimes {\mathcal P}) \stackrel {m}\rightarrow f_*(\omega_{{\mathcal X}/{\mathcal R}^0_g}^{\otimes 2}) \rightarrow 0
\end{equation}
where $m$ is fibrewise the multiplication map and we denote by ${\mathcal I}_2$ the conormal bundle
$ {\mathcal N^{*}_{{\mathcal R}^0_g/{\mathcal A}_{g-1}}}$. Recall that the second fundamental form of
the exact sequence \eqref{I2} is defined as follows
$$II: {\mathcal I}_2 \rightarrow  f_*(\omega_{{\mathcal X}/{\mathcal R}^0_g}^{\otimes 2}) \otimes \Omega^1_{{\mathcal R}^0_g}, \ \
II(s) = m(\nabla(i(s))),$$
where $\nabla$ is the metric connection on
$S^2f_*(\omega_{{\mathcal X}/{\mathcal R}^0_g} \otimes {\mathcal P})=S^2 {\mathcal F}^{Pr}$ defined above.
At the point $(C,A) \in {\mathcal R}^0_g$ the exact sequence \eqref{I2} becomes
$$0 \rightarrow I_2(K_C\otimes A)  \rightarrow S^2H^0(K_C \otimes A) \stackrel {m}\rightarrow H^0(K_C^{\otimes 2})
\rightarrow 0. $$ Hence, if we identify $T_{{\mathcal R}^0_g,b_0}
\stackrel{\cong}\rightarrow H^{1}(T_C)$ via the Kodaira-Spencer
map of the family ${\mathcal X} \stackrel{f} \rightarrow {\mathcal
R}^0_g$, the second fundamental form $II$ at $[C]$ can be seen as
a map $II:I_2(K_C\otimes A)\rightarrow H^0(2K_C)\otimes
H^0(2K_C)$.

\subsection{Gaussian maps}
 Let $Y$ be a smooth complex projective variety and let $\Delta_Y\subset
Y\times Y$ be the diagonal. Let $L$ and $M$ be line bundles on $Y$.
For a non-negative integer $k$, the \emph{k-th Gaussian map}
associated to these data is given by restriction to the diagonal 
\begin{equation}\label{gaussian1}\gamma^k_{L,M}:H^0(Y\times Y,I^k_{\Delta_Y}\otimes
L\boxtimes M )\rightarrow
H^0(Y\times Y,I^k_{\Delta_Y}\otimes
L\boxtimes M  \otimes \OO_{\Delta_Y})\cong
H^0(Y,S^k\Omega_Y^1\otimes L\otimes M).
\end{equation}
The exact sequence
\begin{equation}
\label{Ik} 0 \rightarrow I^{k+1}_{\Delta_Y} \rightarrow
I^k_{\Delta_Y} \rightarrow S^k\Omega^1_Y \rightarrow 0,
\end{equation} twisted by $L\boxtimes M$, shows that the domain of the $k$-th
gaussian map is the kernel of the previous one:
$$\gamma^k_{L,M}:
ker \gamma^{k-1}_{L,M}\rightarrow H^0(Y,S^k\Omega_Y^1\otimes L\otimes
M).$$

In this paper, we will exclusively deal with the second gaussian
map for curves $C$, assuming also that $L =M$.

The map $\gamma^0_L$ is the multiplication map of global
sections
$$H^0(C,L)\otimes
H^0(C,L)\rightarrow H^0(C,L^{\otimes 2})$$
 which obviously
vanishes identically on $\wedge^2 H^0(L)$.
 Consequently, $H^0(C
\times C, I_{\Delta_C}\otimes L\boxtimes L)$ decomposes as
$\wedge^2 H^0(L)\oplus I_2(L)$, where $I_2(L)$ is the kernel of
$S^2H^0(C,L)\rightarrow H^0(C,L^{\otimes 2})$. Since $\gamma^1_L$ vanishes
on symmetric tensors, one writes
$$\gamma^1_L:\wedge^2H^0(L)\rightarrow H^0(K_C\otimes
L^{\otimes 2}).$$
 Again, $H^0(C\times C, I^2_{\Delta_C}\otimes L\boxtimes L)$
 decomposes as  the sum of $I_2(L)$ and the kernel of $\gamma^1_L$. Since $\gamma_L^2$ vanishes identically on skew-symmetric
 tensors, one usually writes
 $$\gamma^2_L:I_2(L)\rightarrow H^0(K_C^{\otimes 2}\otimes L^{\otimes 2})
 $$
Assume now that the line bundle $L$ is $K_C \otimes A$, with $A \in Pic^0(C)[2]$, and denote by
\begin{equation}
\label{muA}
\mu_A:= \gamma^2_{K_C \otimes A}: I_2(K_C \otimes A) \rightarrow H^0(K_C^{\otimes 4})
\end{equation}
 the second gaussian map. It is useful to provide also a local description of it.
Fix a basis $\{\omega_i\}$ of $H^0(K_C\otimes A)$ and write it  in a local coordinate $z$ as $\omega_i = f_i(z)dz \otimes l$, where $l$ is a local generator of the line bundle $A$. For a quadric $Q \in I_2(K_C \otimes A)$ we have $Q = \sum_{i,j} a_{ij} \omega_i \otimes \omega_j $, where $a_{ij} = a_{ji}$ and  $\sum_{i,j} a_{ij} f_i f_j \equiv 0$, hence we have $\sum_{i,j} a_{ij} f'_i f_j \equiv 0$.
The local expression of $\mu_A(Q)$ is
$$\mu_A(Q) =\sum_{i,j} a_{ij} f''_i f_j (dz)^4 = - \sum_{i,j} a_{ij} f'_i f'_j  (dz)^4.$$

The maps $\mu_A$ glue together to give a map of vector bundles on ${\mathcal R}^0_g,$
\begin{equation}
\label{mu}
\mu: {\mathcal I}_2 \rightarrow f_*( (\omega_{{\mathcal X}/{\mathcal R}^0_g}\otimes {\mathcal P})^{\otimes 2}\otimes {\omega_{{\mathcal X}/{\mathcal R}^0_g}}^{\otimes 2})\cong f_*(\omega_{{\mathcal X}/{\mathcal R}^0_g}^{\otimes 4}),
\end{equation}
 where ${\mathcal I}_2 $ is as in \eqref{I2}.

\subsection{The theorem}
In this subsection we show that the second fundamental form $II$
of the Prym map is a lifting of the second Gaussian map $\mu_A$ as
it happens for the second fundamental form of the period map and
the second Gaussian map of the canonical line bundle (see Theorems
2.1 and 4.5 of \cite{cpt}).

To this purpose, recall that given a
holomorphic line bundle $A$ of degree zero on a curve $C$, there
exists a unique (up to constant) hermitian metric $H$ on $A$ and a
unique connection $D_H$ on $A$ which is compatible both with the
holomorphic structure and with the metric and which is flat (see
e.g. \cite{gh}).
If moreover $A^{\otimes 2} = \OO_C$ and we denote by $\pi: \tilde{C} \rightarrow C$ the associated unramified double covering, we can
take an atlas $\{(U_{\alpha}, s_{\alpha})\}$ of $A$ such that the sections $s_{\alpha}$ have values in $\tilde{C}$, hence the cocycle $g_{\alpha,\beta} = s_{\alpha}/s_{\beta}$ has values in $\{\pm1\}$, so it induces   a flat structure on $A$ and a compatible flat hermitian metric on $A$, which is then equal to $H$ up to scalar (see \cite{ko}).

So we can write $D_H = D'_H +
\overline{\partial}$, where $D'_H$ is the $(1,0)$ component. Such
a pair $(A,H)$ is also called a harmonic line bundle and we have
the following properties (see \cite{sim}):

\begin{itemize}
\item The K\"ahler identities.
\item The associated harmonic decomposition
$${\mathcal A}^{\bullet}(A) = {\mathcal H}^{\bullet}(A) \oplus im(D_H) \oplus im(D^*_H)= {\mathcal H}^{\bullet}(A) \oplus
im(\overline{\partial}) \oplus im(\overline{\partial}^*),$$ where
${\mathcal H}(A)$ is the kernel of the laplacian operator $\Delta
= D_H D^*_H + D^*_H D_H =
2(\overline{\partial}\overline{\partial}^*+
\overline{\partial}^*\overline{\partial})$.
 \item The principle of two types
 $$ker(D'_H) \cap ker(\overline{\partial}) \cap (im(D'_H) + im(\overline{\partial})) = im(D'_H \overline{\partial}).$$
\end{itemize}

\begin{TEO}
\label{II}
The diagram
\begin{equation}
 \xymatrix{
 & I_2(K_C \otimes A)\ar[d]^{-\frac{1}{2 \pi i}\mu_A} \ar[r]^{II \ \ } & \ \ S^2 H^0(K_C^{\otimes 2})  \ar[dl]^{\bf{m}}\\
&H^0(K_C^{\otimes 4})
 & }\end{equation}
is commutative.
\end{TEO}
\proof The proof follows the lines of the proof of Th.2.1 of
\cite{cpt}. First of all we take $v\in H^1(C,T_C)$ and we compute
$II(Q)(v)$
 for every $Q\in I_2(K_C \otimes A)$. Using the Kodaira Spencer map $k$ we can assume that
 $v=k(\frac{\partial}{\partial t})$, where $t$ is the local
 coordinate of the unit disc $\Delta=
\{|t|<1\}$ parametrizing a
 one dimensional deformation ${\mathcal X} \stackrel{f}
\rightarrow \Delta$ where $(C,A) = f^{-1}(0)$. Take $Y$ a
${\mathcal C}^{\infty}$ lifting of the holomorphic vector field
$\frac{\partial}{\partial t}$ on $\Delta$, so we have a ${\mathcal
C}^{\infty}$ trivialization $\tau: \Delta \times (C,A) \rightarrow
{\mathcal X}$, $\tau(t,x) := \Phi_{tY}(1)$, where $\Phi_Y(t)$ is
the flow of the vector field $Y$. Then $\theta:=
\overline{\partial} Y_{|(C,A)}$ is  a closed form in
$A^{0,1}(T_C)$ such that $ [\theta]=v\in H^1(C,T_C)$. Denote by
$(C_t, A_t)$ the fibre of $f$ over $t$, where $A_t$ is a
holomorphic line bundle in $Pic^0(C_t)[2]$ endowed with the
flat structure  induced by the double covering $\pi_t: \tilde{C_t} \rightarrow C_t$. We denote by $H_t$ the flat hermitian metric and by
${D_{H_t}} = D'_{H_t}+ \overline{\partial}_t$ the flat Chern connection.

Let $\omega(t)$ be a section of ${\mathcal F}^{Pr}$, hence
$\forall t \in \Delta$, $\omega(t) \in H^0(K_{C_t} \otimes A_t)
\cong H^{1,0}(A_t)$. Denote by $\tau_t: C \rightarrow C_t$  and by by $\sigma_t: \tilde{C} \rightarrow \tilde{C_t}$ the
diffeomorphisms induced by $\tau$, where $\tilde{C}$ and $\tilde{C_t}$ are the unramified double coverings induced by $A$ and by $A_t$. We have the following  commutative diagram:
$$
 \xymatrix{
 &\tilde{C}\ar[d]^{\pi } \ar[r]^{\sigma_t} & \tilde{C_t} \ar[d]^{\pi_t}\\
&C  \ar[r]^{\tau_t}& C_t
 & }$$
so we have an induced map by pullback $\tau_t^{*}: A^1(A_t)
\rightarrow A^1(A)$, and since $\omega_t \in A^{1,0}(A_t)$ is $
D_{H_t}$-closed, then also $\tau_t^{*}(\omega_t)$ is $D_H$-closed
because we have $\tau_t^{*}(D_{H_t}) = D_H$. In fact by the commutativity of the diagram one immediately sees that  the flat structure on $A_t$ given by the covering $\pi_t: \tilde{C_t} \rightarrow C_t$ induces by pullback the flat structure on $A$ given by the covering $\pi: \tilde{C} \rightarrow C$.


So $\tau_t^{*}(\omega_t)$ is
$D_H$-closed, hence it has a power series expansion at $t=0$
$$\tau_t^{*}(\omega_t)= \omega + (\alpha + D_H h)t +o(t),$$
where $\omega := \omega(0)$, $\alpha \in A^1(A)$ is harmonic and
$h$ is a ${\mathcal C}^{\infty}$ section of $A$ (by the harmonic
decomposition for $D_H$). So we have
$\nabla^{GM}_{\frac{\partial}{\partial t}}[\omega(t)]_{t=0} =
[\alpha]$, $\theta \cdot \omega = \alpha^{0,1} +
\overline{\partial}h$, so $k(\frac{\partial}{\partial t}) \cdot
[\omega] = [ \alpha^{0,1}]$, where $  \alpha^{0,1}$ is the $(0,1)$
component of $\alpha$.

Now assume that $\{\omega_i\}_{i=1,...,g-1}$ is a basis of
$H^0(K_C \otimes A) $. Take a quadric $Q \in I_2(K_C \otimes A)$,
$Q = \sum_{i,j} a_{i,j} \omega_i \otimes \omega_j$, with $a_{i,j}
= a_{j,i}$, then $\forall i$ we have
$\nabla^{1,0}_{\frac{\partial}{\partial t}}[\omega_i(t)]_{t=0} =
[{\alpha_i}^{1,0}]$, so if $\tilde{Q}(t) = \sum_{i,j} a_{i,j}(t)
\omega_i(t) \otimes \omega_j(t)$ is a section of ${\mathcal I}_2$
such that $\tilde{Q}(0) = Q$, we have $II(Q)(v) = m (
\nabla_{\frac{\partial}{\partial t}} \tilde{Q}_{|t=0}) =
\sum_{i,j} a'_{i,j}(0) \omega_i \omega_j + 2 \sum_{i,j} a_{i,j}
\alpha^{1,0}_i \omega_j$. Since $\sum_{i,j} a_{i,j}(t) \omega_i(t)
\omega_j(t) \equiv 0$, also its derivative with respect to $t$ at
$t=0$ must be zero, i.e. $2\sum_{i,j} a_{i,j}(\alpha_i + D_H h_i)
\omega_j + \sum_{i,j} a'_{i,j}(0) \omega_i \omega_j \equiv 0$, and
if we take the $(1,0)$ part we have $2\sum_{i,j}
a_{i,j}(\alpha^{1,0}_i + D'_H h_i) \omega_j + \sum_{i,j} a'_{i,j}
\omega_i \omega_j \equiv 0,$ so $II(Q)(v) = -2 \sum_{i,j} a_{i,j}
\omega_j D'_H h_i. $

Now we observe that $ \sum_{i,j} a_{i,j} \omega_j D'_H
h_i=\rho(Q)(v)$ where $\rho$ is the map defined in theorem 4.4. of
\cite{cpt}. So we conclude by  theorem 4.5 of \cite{cpt} that
asserts that $\rho$ is a lifting of $\mu_A$.
 \qed
\begin{COR}
Let $\xi_P \in H^1(T_C)$ be a Schiffer variation at a point $P \in C$. Then we have 
$$\mu_A(Q)(P) = -\frac{1}{2\pi i} II(Q)(\xi_p \odot \xi_P) = -\frac{1}{2\pi i}\xi_P(II(Q)(\xi_P)).$$
\end{COR}
\proof
Recall that given a point $P \in C$, a Schiffer variation $\xi_P \in H^1(T_C)$  is a generator of the image of the coboundary map $H^0(T_C(P)_{|P}) \rightarrow H^1(T_C)$.
Given a local coordinate in a neighborhood of $P$,  under the Dolbeault
isomorphism $H^1(T_X) \cong H^{0,1}(T_X)$, $\xi_P$ is represented by
the form
$\theta_P = \frac{1}{z-z(P)} \overline{\partial}b_P \otimes
\frac{\partial}{\partial z},$ where $b_P$ is a bump function around $P$.
Notice that if we choose $b_P$ to be one in a neighborhood of $P$, $\xi_P$ depends only on
the choice of $z$.  The choice of the local coordinate also allows us to see the evaluation $val_P$ in $P$ as an element of $H^0(K_C^{\otimes 4})^{*}$ and it holds: $m^*(val_P) = \frac{1}{(2 \pi i)^2} \xi_P \odot \xi_P$ (see \cite{cpt} p. 139).
In fact, if $\{\lambda_i \}$ is a basis of $H^0(K_C^{\otimes 2})$, such that locally $\lambda_i = g_i(z) (dz)^2$, we have
$(\xi_P \odot \xi_P)(\lambda_i \odot \lambda_j) = \xi_P(\lambda_i) \xi_P(\lambda_j)$ and by Serre duality
\begin{equation}
\label{conto}
\xi_P(\lambda_i) = \int_{C} g_i(z) dz \wedge \frac{\overline{\partial} b_P}{z-z(P)}    =  -\int_C d \left(\frac{b_Pg_i(z) }{z-z(P)}dz \right) = \int_{\Gamma} \frac{g_i(z)}{z - z(P)} dz=(2\pi i) g_i(P),
\end{equation}
where $\Gamma$ is a small circle around $P$.
Hence  by Theorem \eqref{II} we have $\mu_A(Q)(P) = -\frac{1}{2\pi i} II(Q)(\xi_p \odot \xi_P) = -\frac{1}{2\pi i}\xi_P(II(Q)(\xi_P))$.
\qed

\subsection{Curvature}
In this subsection we would like to give an explicit formula for the
holomorphic sectional curvature of the Siegel metric on ${\mathcal
R}_g$ along the tangent directions given by the Schiffer
variations. This formula is analogous to the formula of the
holomorphic sectional curvature of the Siegel metric on ${\mathcal
M}_g$ induced by the Period map, given in Cor.3.8 \cite{cf2}.\footnote{Notice that between the formula of Cor 3.8 of \cite{cf2} and the formula given in Prop. \ref{curv} there is a difference by a factor 4, due to a small mistake in  Cor 3.8 of \cite{cf2} where we identified $\rho$ with $II$, while $II = -2 \rho$ (\cite{cpt} p.136).  }

Assume that $\{Q_i\}$ is an orthonormal basis of $I_2(K_C\otimes
A)$, $\{\omega_i\}$ an orthonormal basis of of $H^0(K_C\otimes A)$
and choose a local coordinate $z$ at $P$ and a local generator $a$
of $A$ such that locally $\omega_i=f_i(z)dz\otimes a$.

\begin{PROP}
\label{curv}
The holomorphic sectional curvature $H$ of $T_{{\mathcal R}^0_g}$
at $[C,A]\in {\mathcal R}^0_g$ computed at the tangent vector
$\xi_P$ given by a Schiffer variation in $P$ is given by:
$$H(\xi_P)=-1-\frac{1}{16\alpha_{P}^4 \pi^2}\sum_i|\mu_A(Q_i)(P)|^2$$
where   $\alpha_{P}=\sum_i|f_i(P)|^2$.
\end{PROP}

\proof Also the proof follows the lines of  \cite{cf2}. We start
recalling that
$$H(\xi_P) =\frac{\langle R(\xi_P), \xi_P\rangle(\xi_P, \overline{\xi_P}) }{{\langle \xi_P, \xi_P \rangle}^2},$$
and by the Gauss formula we have $\langle R(\xi_P),
\xi_{P}\rangle(\xi_P, \overline{\xi_P}) = \langle
\tilde{R}(\xi_P), \xi_{P}\rangle (\xi_P,
\overline{\xi_P})  - \langle\sigma(\xi_P),
\sigma(\xi_{P})\rangle(\xi_P, \overline{\xi_P}),$ where $R$ is the curvature of the Siegel metric, $\tilde{R}$ is the curvature of ${\mathcal A}_{g-1}$ and $\sigma$ is the second fundamental form of ${\mathcal R}_g^0$ in ${\mathcal A}_{g-1}$ (see the tangent bundle exact sequence \eqref{tangent}).

First of all observe that the holomorphic sectional curvature of ${\mathcal A}_{g-1}$ along the Schiffer variations is equal to $-1$ (see the argument below corollary 3.8 of \cite{cf2}). In fact a Schiffer variation $\xi_P$, seen as a symmetric homomorphism $H^0(K_C \otimes A) \rightarrow  H^0(K_C \otimes A)^* \cong H^1(A)$ has rank 1, since its kernel is $H^0(K_C \otimes A(-P))$.
To compute $\langle\sigma(\xi_P), \sigma(\xi_{P})\rangle(\xi_P, \overline{\xi_P})$, recall that $\sigma(\xi_P)(Q) = II(Q)(\xi_P)$, for all $Q \in I_2(K_C\otimes A)$, hence $\sigma(\xi_P) = \sum_i II(Q_i)(\xi_P) \otimes Q^*_i$, where $\{Q_i\}$ is an orthonormal basis of $I_2(K_C\otimes A)$. Hence $\langle\sigma(\xi_P),
\sigma(\xi_{P})\rangle =  \sum_i \langle II(Q_i)(\xi_P) , II(Q_i)(\xi_P) \rangle$.  Now recall that the Schiffer variations at $3g-3$ general points of $C$ give a basis of $H^1(T_C)$, so we can write $II(Q_i)(\xi_P) = \sum_S \xi_S(II(Q_i)(\xi_P)) \xi_S^*$, therefore $ \langle\sigma(\xi_P), \sigma(\xi_{P})\rangle(\xi_P, \overline{\xi_P}) = \sum_i (\xi_P(II(Q_i)(\xi_P))\overline{( \xi_P(II(Q_i)(\xi_P))} .$ By \ref{schiffer} we have $\sum_i (\xi_P(II(Q_i)(\xi_P))\overline{( \xi_P(II(Q_i)(\xi_P) )}=  4 \pi^2 \sum_i|\mu_A(Q_i)(P)|^2$. Now it remains to show that $\langle \xi_P, \xi_P \rangle = 8 \pi^2 \alpha_P^2$ as in Lemma 2.2 of \cite{cf2}. To do this we write $\xi_P$ as an element of $S^2(H^0(K_C \otimes A)^*)$ as $\xi_P = \sum_{i,j} \xi_P(\omega_i)(\omega_j)(\omega_j^* \odot \omega_i^*)$, where $\{\omega_i^*\}$ is the dual basis of the orthonormal basis $\{\omega_i\}$. One computes $\xi_P(\omega_i)(\omega_j) = \xi_P(\omega_i \omega_j) = 2 \pi i f_i(P) f_j(P)$ by \eqref{conto}, then the proof follows exactly as in lemma 2.2 of \cite{cf2}.
\qed

\section{Prym-canonical binary curves}
\subsection{Strategy of the proof of surjectivity}

The rest of the paper is devoted to the proof of the surjectivity
of the 2nd gaussian map $\mu_A$ for the general point $[C,A] \in
{\mathcal R}_g$. We will do it by degeneration to binary curves
following the method used in \cite{ccm} for the second gaussian
map of the canonical line bundle. We recall that ${\mathcal R}_g$
admits a suitable compactification $\overline{{\mathcal R}}_g$,
which is isomorphic to the coarse moduli space of the stack
${\bf{R}}_g$ of Beauville admissible double covers (\cite{b},
\cite{acv}) and to the coarse moduli space of the stack of Prym
curves (\cite{bcf}).

Consider the partial compactification $\tilde{\mathcal R}_g$ of ${\mathcal R}_g$ introduced in \cite{fl}. Denote by $f: {\mathcal X} \rightarrow \tilde{{\bf{R}}}_g$ the universal family and by ${\mathcal P} \in Pic({\mathcal X})$ the corresponding Prym bundle as in \cite{fl} 1.1.
The map of vector bundles over ${\mathcal R}^0_g$, $\mu: {\mathcal I}_2 \rightarrow f_*( (\omega_{{\mathcal X}/{\mathcal R}^0_g}\otimes {\mathcal P})^{\otimes 2}\otimes {\omega_{{\mathcal X}/{\mathcal R}^0_g}}^{\otimes 2})\cong f_*(\omega_{{\mathcal X}/{\mathcal R}^0_g}^{\otimes 4})$ defined in \eqref{mu}, extends to a map
\begin{equation}
\label{mutilde}
\tilde{\mu}: \tilde{{\mathcal I}}_2 \rightarrow f_*( (\omega_f\otimes {\mathcal P})^{\otimes 2}\otimes S^2({\Omega^1_f}))\cong f_*(\omega_f^{\otimes 4}\otimes {\mathcal P}^{\otimes 2}\otimes {\mathcal I}_Z^{\otimes 2}),
\end{equation}
where $\tilde{{\mathcal I}}_2 $ is the kernel of the
multiplication map $ S^2f_*( \omega_f \otimes {\mathcal P})
\rightarrow f_{*}(\omega_f^{\otimes 2} \otimes {\mathcal
P}^{\otimes 2}),$ and $Z$ is the locus of nodes of fibres of $f$, so $\Omega^1_f \cong \omega_f
\otimes {\mathcal I}_Z$.

If $[C,A] \in $ is a point in $\tilde{\mathcal R}_g$, the local expression of
\begin{equation}
\label{gauss}
\mu_A: I_2(\omega_C \otimes A) \rightarrow H^0((\omega_C \otimes A)^{\otimes 2} \otimes S^2({\Omega^1_C}))
\end{equation}
 is as follows.  Let  $\{\omega_i\}$ be a basis of $H^0(\omega_C\otimes A)$ and write it in a local coordinate as $\omega_i = f_i(z)\xi \otimes l$, where $\xi$ and $l$ are local generators of the line bundles $\omega_C$, respectively $A$.
For a quadric $Q = \sum_{i,j} a_{ij} \omega_i \otimes \omega_j \in I_2(\omega_C \otimes A)$, $\mu_A(Q)$ is locally defined as
\begin{equation}
\label{mulocal}
\mu_A(Q) = - \sum_{i,j} a_{ij} (df_i)  (df_j)  \xi^{\otimes2} \otimes l^{\otimes 2}.
\end{equation}
To prove by semicontinuity the surjectivity of $\mu_A$ for the general point in ${\mathcal R}_g$ in the following we will exhibit a Prym-canonical binary curve $(C,A)$ for which $\mu_A$ is surjective.

\subsection{Construction of Prym-canonical binary curves}
Recall that a binary curve of genus $g$ is a stable curve
consisting of two rational components $C_j$, $j=1,2$ meeting
transversally at $g+1$ points. Moreover one can check that if
$A\in Pic^0(C)$ then $H^0(C,\omega_C \otimes A)$ has dimension $g-1$
and  the restriction of $\omega_C \otimes A$ to  the component $C_j$ is
$K_{C_j}(D_j)$ where $D_j$ is the divisor of nodes on $C_j$. Since
$K_{C_j}(D_j)\cong \OO_{{\proj}^{1}}(g-1)$ we observe that the
components are embedded by a linear subsystem of
$\OO_{{\proj}^{1}}(g-1)$, hence they are projections from a point
of rational normal curves in ${\proj}^{g-1}$. Viceversa, let us
take 2 rational curves embedded in  ${\proj}^{g-2}$ by non
complete linear systems of degree $g-1$ intersecting transversally
at $g+1$ points. Then their union $C$ is a binary curve of genus
$g$ embedded either by a linear subsystem of $\omega_C$ or by a
complete linear system $|\omega_C \otimes A|$, where $A\in Pic^0(C)$
is nontrivial (see e.g. \cite{capo}, Lemma 10). In this section we will construct a
binary curve $C$ embedded in ${\proj}^{g-2}$ by  a linear system
$|\omega_C \otimes A|$ with $A^{\otimes 2}\cong \OO_C$, and $A$ is non
trivial.

Assuming that the first $g-1$ nodes, $P_1,...,P_{g-1}$  are in
general position, up to projective transformations we will take
$P_i=(0,...,0,1,0,...0)$ with 1 at the $i$-th place.  Then we can
assume that $C_j$ is the image of the map
\begin{equation}\label{phi}
\begin{gathered}\phi_j:{\proj}^1 \rightarrow {\proj}^{g-2}, \
j=1,2\\
\phi_j(t,u):= [M_j(t,u)\frac{(\de_{1,j}t-c_{1,j}u)}{(t-a_{1,j}u)},
..., M_j(t,u)\frac{(\de_{g-1,j}t-c_{g-1,j}u)}{(t-a_{g-1,j}u)}]
\end{gathered}
\end{equation}
with $M_j(t,u):= \prod_{r=1}^{g-1} (t-a_{r,j}u)$, $j=1,2$ and
$\phi_j([a_{l,j},1]) = P_l$, $l=1,...,g-1$.

We will also impose  that the remaining two nodes
$P_g:=[t_1,...,t_{g-1}]$ and $P_{g+1}:=[s_1,...,s_{g-1}]$ are the
images of $[0,1]$ and $[1,0]$ through the maps $\phi_j$, $j=1,2$.
This is equivalent to 
$$\ c_{i,j}= \frac{d_j t_ia_{i,j}}{A_j}, \
\de_{i,j}= \mu_js_i, $$ where $\mu_j, d_j$ are non zero scalars and
$A_j= \prod_{k=1}^{g-1} a_{k,j}$, $j=1,2$.

\begin{LEM}
\label{lemmadeltac} Let us choose $s_i=1$ for $i
=1,...,[\frac{g}{2}]$, $s_i =0$, for $i =
[\frac{g}{2}]+1,...,g-1$, while $t_i=0$ for $i
=1,...,[\frac{g}{2}]$, $t_i =1$, for $i =
[\frac{g}{2}]+1,...,g-1$, $\mu_1 = \mu_2 =: \mu$, $d_1=-\frac{d_2
A_1}{A_2}$. Then, for a general choice of $a_{i,j}$'s, $C=C_1\cup
C_2$ is a binary curve embedded in ${\proj}^{g-2}$ by  a linear
system $|\omega_C \otimes A|$ with $A^{\otimes 2}\cong \OO_C$ and $A$ nontrivial.
\end{LEM}

\proof

One can easily check that if we choose the elements $a_{k,j}$
general ($j=1,2$, $k=1,...,g-1$), $C_j$ are smooth rational curves
and $C$ has exactly $g+1$ nodes at the points $P_k$, $k
=1,...,g+1$, and no other singularity. Then, by the above
discussion we know that $C$ is a binary curve  embedded in
${\proj}^{g-2}$ by a linear system of $\omega_C \otimes A$, with
$deg(A)=0$.
 We will now show that $A$ is a 2-torsion non trivial element in
$Pic^0(C)$. In fact, recall that $Pic^0(C) \cong {{\Co}^*}^g$ and
if we denote by $\alpha:N \rightarrow C$ the normalization map, we
have an exact sequence
\begin{equation}
\label{KA} 0 \rightarrow (\omega_C\otimes A) \rightarrow
\alpha_*(\alpha^*(\omega_C\otimes A)) \rightarrow
\oplus_{i=1}^{g+1}{\Co}_{P_i}  \rightarrow 0.
\end{equation}
If we set $\{q_i,r_i \} = \alpha^{-1}(P_i)$, with $q_i \in C_1$,
$r_i \in C_2$, $i=1,...,g+1$, $D_1:=\sum_{i=1}^{g+1}q_i$,
$D_2:=\sum_{i=1}^{g+1}r_i$ we have $\alpha^*(\omega_C\otimes A) =
K_N(D_1+D_2)$. So if we take the long exact sequence in cohomology
associated to (\ref{KA}), we have
\begin{equation}
0 \rightarrow H^0(\omega_C\otimes A) \rightarrow H^0(K_{C_1}(D_1))
\oplus H^0(K_{C_2}(D_2)) \stackrel{e}\rightarrow {\Co}^{g+1}
\rightarrow 0.
\end{equation}
Clearly $ H^0(K_{C_1}(D_1)) \cong H^0({\mathcal
O}_{{\proj}^1}(g-1)) \cong {\Co}^g \cong H^0(K_{C_2}(D_2))$.
Recall that the line bundle $A$ corresponds to an element in $
{{\Co}^*}^g$ as follows. Consider the natural isomorphisms $f_i:
({\alpha^*(A)_{|C_1}})_{q_i} \rightarrow
({\alpha^*(A)_{|C_2}})_{r_i}$, and choose local trivializations $
({\alpha^*(A)_{|C_1}})_{q_i} \cong {\Co}$, $
({\alpha^*(A)_{|C_2}})_{r_i} \cong {\Co}$, $\forall i =1,...,g+1$.
Then $f_i$ is given by multiplication by an element $h_i \in
{\Co}^*$, $\forall i = 1,...,g+1$, and we associate to $A$ the
element $(h_1,...,h_{g+1})$ modulo the diagonal action of
${\Co}^*$.

Notice that if $\sigma \in H^0(\omega_C\otimes A)$ and
$\alpha^*(\sigma) = (\sigma_1, \sigma_2) \in   H^0(K_{C_1}(D_1))
\oplus H^0(K_{C_2}(D_2))$, then  we have
$$Res_{q_i}(\sigma_1) - h_i Res_{r_i}(\sigma_2) =0, \ \forall i=1,...,g+1.$$
We claim that with our assumptions the line bundle $A$ corresponds
to the element $[(h_1,...,h_{g+1})] \in {{\Co}^*}^{g+1}/{\Co}^*$,
where $h_i=1$, for $i< [\frac{g}{2}]+1$, $h_i = -1$, for $i
=[\frac{g}{2}]+1,...,g-1$, $h_g=-1$, $h_{g+1} = 1$, so $A$ is of
2-torsion. In fact, consider the hyperplane $x_i =0$,
$i=1,...,g-1$ in ${\proj}^{g-2}$, and set $\sigma_{i,1} :=
\phi_{1}^*(x_i) \in H^0({\mathcal O}_{{\proj}^1}(g-1)) \cong
H^0(K_{C_1}(D_1))$, $\sigma_{i,2} := \phi_{2}^*(x_i) \in
H^0({\mathcal O}_{{\proj}^1}(g-1)) \cong H^0(K_{C_2}(D_2))$. We
have
$$\sigma_{i,j} = \frac{(\de_{i,j}t -c_{i,j})}{(t - a_{i,j})t} dt, \ j=1,2.$$
Notice that, with our assumptions, we have
\begin{equation}
\label{deltac}
\begin{gathered}
\de_{i,1} = \de_{i,2} = \mu, \ i < [\frac{g}{2}]+1,\ \de_{i,1} = \de_{i,2} = 0,  \  i \geq [\frac{g}{2}]+1,\\
c_{i,1} = c_{i,2} =0, \ i< [\frac{g}{2}]+1,\\
c_{i,1} = \frac{{d_1} a_{i,1}}{A_1}= -\frac{{d_2}a_{i,1}}{A_2}, \
c_{i,2} = \frac{{d_2} a_{i,2}}{A_2}, \ i \geq [\frac{g}{2}]+1,
\end{gathered}
\end{equation}
and for simplicity we shall choose $d_2 =\mu = 1$, so the only parameters are the $a_{i,j}$'s.
Hence, for $j =1,2$, we have $Res_{q_i}(\sigma_{i,1}) = \de_{i,1}
- c_{i,1}/a_{i,1} = \mu = \de_{i,2} - c_{i,2}/a_{i,2}=
Res_{r_i}(\sigma_{i,2})$ for $i  < [\frac{g}{2}]+1$, so
$$h_i = \frac{Res_{q_i}(\sigma_{i,1})}{Res_{r_i}(\sigma_{i,2}) }= 1, \ for \ i < [\frac{g}{2}]+1,$$
$$h_i =  \frac{Res_{q_i}(\sigma_{i,1})}{Res_{r_i}(\sigma_{i,2}) }= \frac{\frac{c_{i,1}}{a_{i,1}}}{ \frac{c_{i,2}}{a_{i,2}}}= -1, \ for \ i = [\frac{g}{2}]+1,..., g-1,$$
$$h_{g}= \frac{Res_{q_{g}}(\sigma_{g-1,1})}{Res_{r_{g}}(\sigma_{g-1,2}) }= \frac{\frac{c_{g-1,1}}{a_{g-1,1}}}{\frac{c_{g-1,2}}{a_{g-1,2}}}= -1,$$
$$h_{g+1} = \frac{Res_{q_{g+1}}(\sigma_{1,1})}{Res_{r_{g+1}}(\sigma_{1,2}) } = \frac{Res_0(    \frac{(\de_{1,1}-u c_{1,1})} {(1 - \de_{1,1}u)} (-\frac{1}{u}) du)}{ Res_0(    \frac{(\de_{1,2}-u c_{1,2})} {(1 - \de_{1,2}u)} (-\frac{1}{u}) du)}= \frac{\de_{1,1}}{\de_{1,2}} = 1.$$
\qed

\section{Quadrics}
In this section we explicitly describe the ideal
$I_2(C) := I_2(\omega_C \otimes A)$ of the quadrics containing the Prym-canonical binary curve $C$ embedded in ${\proj}^{g-2}$ by $\omega_C \otimes A$ as in the previous section for a general choice of the $a_{i,j}$'s. Similarly as in Proposition 7 of \cite{ccm},  the ideal $I_2(C)$ is described as the space of solutions of the linear system given in Proposition \eqref{Z} which has maximal rank $2g-2$, so the curve is quadratically normal.

Observe that, since the
curves $C_1$ and $C_2$ pass through the coordinate points, the
equation of a quadric $Q \subset {\proj}^{g-2}$ containing $C_k$
has the form
\begin{equation}
\label{Qua} \sum_{1\leq i<j\leq g-1} s_{ij} x_i x_j =0.
\end{equation}

In the next lemma we give a set of generators of $I_2(C_k)$ of the
above form.

\begin{LEM}
Set
\begin{equation}
Q_{n,k}:= \sum_{1 \leq i<j\leq g-1} \tilde{q}_{g-1-n,k;i,j} \cdot
s_{ij}, \ n=0,...,g-1,\ k=1,2
\end{equation}
with
\begin{equation}
\begin{gathered}
\tilde{q}_{0,k;i,j} := q_{0,k;i,j} \de_{i,k} \de_{j,k} \\
\tilde{q}_{1,k;i,j} := q_{1,k;i,j} \de_{i,k} \de_{j,k}- q_{0,k;i,j}(\de_{i,k}c_{j,k} + c_{i,k} \de_{j,k}) \\
\tilde{q}_{r,k;i,j} := q_{r,k;i,j} \de_{i,k} \de_{j,k}-
q_{r-1,k;i,j}(\de_{i,k}c_{j,k} + c_{i,k} \de_{j,k}) +
q_{r-2,k;i,j} c_{i,k} c_{j,k},\ (r\geq 2),
\end{gathered}
\end{equation}
where
\begin{equation}
q_{h,k;i,j}:= \sum_{m=0}^h a_{i,k}^m a_{j,k}^{h-m}
\end{equation}
and $\de_{i,k}$, $c_{i,k}$, $a_{i,k}$ are as in \eqref{phi}.

Then the quadrics in $I_2(C_k)$ ($k=1,2$) are the solutions of the
linear system:
\begin{equation}\label{Qsistem} Q_{n,k}(s_{ij}) = 0, \ \
n=0,...,g-1.
\end{equation}
\end{LEM}

\proof The quadrics of the form \eqref{Qua} containing $C_k$ are
the quadrics which satisfy  the equations:
\begin{equation}
\label{quadric} P_k(t)= \sum_{1\leq i<j\leq g-1}
M_k(t,1)\frac{(\de_{i,k}t - c_{i,k})(\de_{j,k}t -
c_{j,k})}{(t-a_{i,k})(t-a_{j,k})} s_{ij} = \sum_{n=0}^{g-1}
P_{n,k}(s_{ij}) t^n \equiv 0,
\end{equation}
$k=1,2$, where the coefficients $P_{n,k}(s_{ij})$ of the
polynomial $P_k(t)$ are linear in the $s_{ij}$'s. We will show
that the linear system $P_{n,k}(s_{ij})=0$ is equivalent to the
system \eqref{Qsistem}.

By expanding the product $M_k(t,1)$ one sees that the coefficients
$p_{h,k;i,j}$ of $s_{ij}$ in $P_{g-1-h, k}$ are
\begin{equation}
p_{0,k;i,j}= \de_{i,k}\de_{j,k}, \  p_{1,k;i,j}= - \sum_{i_1 \neq i,j} a_{i_1,k} \de_{i,k}\de_{j,k} -(\de_{i,k}c_{j,k} + c_{i,k} \de_{j,k})
\end{equation}

\begin{equation}
\label{pij}
\begin{gathered}
p_{h,k;i,j}= (-1)^h (\sum_{\substack{1\leq i_1 < i_2 < \cdots < i_h \leq g-1\\ \text{all} \ne i,j}} a_{i_1,k}\cdot \cdot \cdot a_{i_h,k} )\de_{i,k}\de_{j,k}  \\
+ (-1)^h  (\sum_{\substack{1\leq i_1 < i_2 < \cdots < i_{h-1} \leq g-1\\ \text{all} \ne i,j}} a_{i_1,k}\cdot \cdot \cdot a_{i_{h-1},k} )(\de_{i,k}c_{j,k} + c_{i,k} \de_{j,k}) \\
+ (-1)^h (\sum_{\substack{1\leq i_1 < i_2 < \cdots < i_{h-2} \leq g-1\\
\text{all} \ne i,j}}a_{i_1,k}\cdot \cdot \cdot a_{i_{h-2},k}
)c_{i,k} c_{j,k},
\end{gathered}
\end{equation}
for $ 2 \leq  h \leq g-1$.

Set
\begin{equation}
\gamma_{0,k} = 1,  \ \gamma_{h,k} = (-1)^h \sum_{\substack{1\leq i_1 < i_2 < \cdots < i_h\leq g-1\\ \text{all} \ne i,j}}a_{i_1,k}\cdot \cdot \cdot a_{i_h,k} ,
\end{equation}

Then we have (cf. \cite{ccm} (12))
\begin{equation}
\label{formula12} (-1)^h \sum_{\substack{1\leq i_1 < i_2 < \cdots
< i_h \leq g-1\\ \text{all} \ne i,j}} a_{i_1,k}\cdot \cdot \cdot
a_{i_h,k} = \sum_{l=0}^h \gamma_{l,k} q_{h-l,k;i,j}
\end{equation}
for $h = 0,...,g -1$.

So, by (\ref{formula12}), formula (\ref{pij}) becomes

\begin{equation}
\label{pijneu}
\begin{gathered}
p_{h,k;i,j}=\sum_{l=0}^{h-2}\gamma_{l,k}(q_{h-l,k;i,j} \de_{i,k}\de_{j,k} -q_{h-1-l,k;i,j} (\de_{i,k}c_{j,k}+c_{i,k}\de_{j,k})+q_{h-2-l,k;i,j} c_{i,k}c_{j,k}) +\\
\gamma_{h-1,k}(q_{1,k;i,j} \de_{i,k}\de_{j,k} -q_{0,k;i,j}(\de_{i,k}c_{j,k}+c_{i,k}\de_{j,k})) +
 \gamma_{h,k} \de_{i,k}\de_{j,k}= \sum_{l=0}^h \gamma_{l,k} \tilde{q}_{h-l,k;i,j} \end{gathered}
\end{equation}
So, we have
\begin{equation}
P_{n,k}= \sum_{m=0}^{g-1-n} \gamma_{m,k} Q_{n+m,k}
\end{equation}
and one immediately checks that the linear systems
$P_{n,k}(s_{ij}) =0$ and \eqref{Qsistem} are equivalent. \qed

\begin{PROP}\label{quad1}
Let $g \geq 6$. For a general choice of $a_{i,k}$, $k=1,2$, $i =
1,...,g-1$ and with conditions \eqref{deltac} on
$c_{i,k},\de_{i,k}$, the linear system (\ref{Qsistem}) has maximal
rank $g$.
\end{PROP}
\proof

Consider the matrix
$$M(a_{1,k},...,a_{g-1,k}):= (\tilde{q}_{h,k;i,j})_{0\leq h\leq g-1, 1\leq i<j\leq g-1}$$
of size $g \times \frac{(g-1)(g-2)}{2}$. We will show that  the
minor $B_g$ determined by the columns with indexes $(i,j)=
(1,2),...,(1,g-1),(2,[\frac{g}{2}] + 1),(g-2,g-1)$ is non zero.

Notice that
\begin{equation}
\begin{gathered}
\tilde{q}_{0,k;1,j} = \mu^2, \   \ \tilde{q}_{1,k;1,j} = \mu^2 q_{1,k;1,j}, \  \tilde{q}_{h,k;1,j} = \mu^2 q_{h,k;1,j}, \  \ j < [\frac{g}{2}] + 1, \ h \geq 2\\
\tilde{q}_{0,k;1,j} = 0, \  \  \tilde{q}_{1,k;1,j} = -\mu c_{j,k}, \  \tilde{q}_{h,k;1,j} = -\mu c_{j,k} q_{h-1,k;i,j}, \  \ j \geq [\frac{g}{2}] + 1, \ h \geq 2\\
\tilde{q}_{0,k;2,[\frac{g}{2}] + 1} = \tilde{q}_{0,k;g-2,g-1}= 0,  \ \tilde{q}_{1,k;2,[\frac{g}{2}] + 1} = -\mu c_{[\frac{g}{2}] + 1,k} , \ \tilde{q}_{1,k;g-2,g-1} = 0 ,\\
\tilde{q}_{h,k;2,[\frac{g}{2}] + 1} = -\mu c_{[\frac{g}{2}] + 1,k}q_{h-1,k;2,[\frac{g}{2}] + 1}, \ h \geq 2\\
\tilde{q}_{h,k;g-2,g-1} = q_{h-2,k;g-2,g-1} c_{g-2,k} c_{g-1,k}, \ h \geq 2 \\
\end{gathered}
\end{equation}
So, dividing the first $[\frac{g}{2}] $ columns by $\mu^2$, the
column indexed by $(2,[\frac{g}{2}] + 1)$ by $-\mu
c_{[\frac{g}{2}] + 1,k}$,the last column by $c_{g-2,k} c_{g-1,k}$
and  all the other columns indexed by $(1,j)$ with $[\frac{g}{2}]
+ 1 \leq j \leq g-1$ by $- \mu c_{j,k}$, we see that $B_{g}$ is a
non zero multiple of the determinant $d$ of the following matrix:

\begin{equation}\label{Mtilde} {\small \left( \begin{array}{cccccccc}
1&..&1&0&..&0&0&0 \\
q_{1,k;1,2} &..&q_{1,k;1,[\frac{g}{2}]}&1&..&1&1&0\\
q_{2,k;1,2} &..&q_{2,k;1,[\frac{g}{2}]}&q_{1,k;1,[\frac{g}{2}]+ 1}&..&q_{1,k;1,g-1}&q_{1,k;2,[\frac{g}{2}] + 1}&1\\
q_{3,k;1,2} &..&q_{3,k;1,[\frac{g}{2}]}&q_{2,k;1,[\frac{g}{2}]+ 1}&..&q_{2,k;1,g-1}&q_{2,k;2,[\frac{g}{2}] + 1}&
q_{1,k;g-2,g-1}\\
.&..&.&.&..&.&.&. \\
.&..&.&.&..&.&.&. \\
.&..&.&.&..&.&.&. \\
q_{g-1,k;1,2} &..&q_{g-1,k;1,[\frac{g}{2}]}&q_{g-2,k;1,[\frac{g}{2}]+ 1}&..&q_{g-2,k;1,g-1}&q_{g-2,k;2,[\frac{g}{2}]+ 1}&q_{g-3,k;g-2,g-1}\\
\end{array} \right)}
\end{equation}

One can inductively compute the determinant $d$ up to sign,
\begin{equation}
\label{det1} d= V(a_{3,k},...,a_{g-1,k})\cdot
\prod_{\substack{r=1\\ r \ne 2,[\frac{g}{2}]+1}}^{g-1}
(a_{r,k}-a_{2,k})\cdot \prod_{s=3}^{[\frac{g}{2}]} a_{s,k} \cdot
\prod_{\substack{j=1}}^{g-3} a_{j,k},
\end{equation}
where $V(a_{3,k},...,a_{g-1,k})$ is the Vandermonde determinant in
the variables $a_{3,k},...,a_{g-1,k}$. To do this one can perform
column and row operations.  \footnote{ Substitute column
$(2,[\frac{g}{2}])$ with $(2,[\frac{g}{2}])-(1,[\frac{g}{2}])$
then divide it by $a_{2,k}-a_{1,k}$; substitute any row by itself
minus $a_{1,k}$ times the preceding row;
 substitute each column from $(1,3)$ to $(1,[\frac{g}{2}])$
by itself minus the first column, eliminate the first row and
column and divide the column $(1,i)$ ($j=3,...,[\frac{g}{2}]$) by
$a_{i,k}-a_{2,k}$; substitute each row by itself minus
$a_{[\frac{g}{2}]+1,k}$ times the preceding row,  eliminate the
first row and the $(1,[\frac{g}{2}]+1)$-column and divide the
column $(1,j)$ with $j=[\frac{g}{2}]+2,...,g-1$ by
$a_{j,k}-a_{[\frac{g}{2}]+1,k}$; substitute any row by itself
minus $a_{2,k}$ times the preceding row, eliminate the first row
and the $(2,[\frac{g}{2}]+1)$-column and divide the column $(1,j)$
with $j=3,...,[\frac{g}{2}]$ by $a_{j,k}-a_{[\frac{g}{2}]+1,k}$.}

So we reduce to a $(g-3)\times (g-3)$ matrix whose columns except
the last one are the columns of the Vandermonde matrix in the
variables $a_{3,k},...,a_{[\frac{g}{2}],k},$
$a_{[\frac{g}{2}]+2,k},...,a_{g-1,k}$. Hence repeating
recursively standard row and columns
operations,\footnote{substitute each column except the first one
and the last one by itself minus the first column, substitute the
last one by itself minus the first column multiplied by the first
coefficient of the last column, eliminate the first row and column
and divide all the columns except the last one by
$a_{j,k}-a_{3,k}$ and repeat.} we obtain formula (\ref{det1}). \qed

In the following proposition we give an explicit description of
the ideal $I_2(C)$ of quadrics containing $C=C_1\cup C_2$ and we
prove that $C$ is quadratically normal.
\begin{PROP}\label{Z}
Let $g \geq 6$. For a general choice of $a_{i,k}$, $k=1,2$, $i =
1,...,g-1$ and with conditions \eqref{deltac} on
$c_{i,k},\de_{i,k}$, the linear system
\begin{equation}\label{qqij} Q_{0,1}(s_{ij})=...=Q_{g-1,1}(s_{ij})=Q_{1,2}(s_{ij})=...=Q_{g-2,2}(s_{ij})=0,\end{equation}
has maximal rank $2g-2$.
\end{PROP}
\proof
Since we want to prove the statement for generic $a_{i,j}$, it suffices to show it for the following choice of $a_{i,j}$, $j=1,2$, $i=1,...,g-1$:
\begin{equation}\label{scelta}
a_{i,1} := i\cdot a,  i=1,...,g-1; \ a_{1,2}:= 1, \ a_{r,2}:=r+1, \ r=1,...,g-1,
\end{equation}
where $a\neq 1$ is a non zero constant. Consider the matrix
$Z:=Z(a_{i,j})$ of size $(2g-2) \times \binom{g-1}{2}$ obtained by
concatenating  vertically $M(a_{1,1},...,a_{g-1,1})=
(\tilde{q}_{h,1;i,j})_{0\leq h\leq g-1, 1\leq i<j\leq g-1}$, $N(a_{1,2},...,a_{g-1,2})$
$= (\tilde{q}_{h,2;i,j})_{1\leq h\leq g-2,
1\leq i<j\leq g-1}$. Let us set $k:=[\frac{g}{2}]$ and consider
the submatrix $Z_1$ of $Z$ formed by the columns of $Z$ indexed by
$(1,2),...,(1,g-1),(2,k+1),(g-2,g-1),(2,3),...,(2,k),(2,
k+2),...,(2,g-1),(k, k+1), (k, g-1)$. We will prove that $Z_1$
has maximal rank $2g-2$. Note that the submatrix  given by the
first $g$ rows and columns is the matrix \eqref{Mtilde} of
Proposition \ref{quad1} which is proved to be non singular. So
doing operations on the columns we can assume that $Z_1$ is a
matrix whose submatrix given by the last $g-2$ columns ad the
first $g$ rows is zero. Hence we just need to prove that the
submatrix $A$ given by the last $g-2$ rows and columns has maximal
rank. If we denote by $v_i$ the column indexed by $(1,i)$, $i
=1,...,g-1$, by $w_i$ the column indexed by $(2,i)$, $i
=1,...,g-1$, by $w$ the column indexed by $(k,k+1)$, by $\zeta$
the column indexed by $(k,g-1)$,  the operations that we do on the
columns of $Z_1$ are the following:

- for $i=3,...,k$, substitute the column $w_i$ with the vector
$w_i+\frac{1-i}{i-2}v_i + \frac{1}{i-2}v_1$.

- for $i=k+2,...,g-1$, substitute the column $w_i$ with the vector

$w_i+\frac{k\cdot c_{i,1}}{c_{k+1,1}(i-2)}v_{k+1} -\frac{c_{i,1}(k-1)}{c_{k+1,1}(i-2)}w_{k+1}- \frac{i-1}{i-2} v_{i}.$

- substitute the column $w$  with the vector

$w+\frac{(k-1)\cdot c_{k+1,1}}{ka}v_{1} -\frac{c_{k+1,1}(k-1)}{ka}v_{k}- (2-k)v_{k+1}-\frac{2(1-k)^2}{k}w_{k+1}$.

- substitute the column $\zeta$  with the vector

$\zeta+\frac{(k-1)\cdot c_{g-1,1}}{ka(g-k-1)}v_{1} -\frac{c_{g-1,1}(k-1)}{ka(g-k-1)}v_{k}+ \frac{2(k-1)c_{g-1,1}}{(g-k-1)c_{k+1,1}}v_{k+1}-\frac{2(1-k)^2c_{g-1,1}}{(g-k-1)kc_{k+1,1}}w_{k+1}- \frac{g-2}{g-k-1}v_{g-1}$.

To prove that the matrix $A$ is of maximal rank $g-2$, we argue
as follows.  First of all one can easily check (with the same
procedure as in  \ref{quad1}) that the submatrix $C$ of $Z_1$
formed by the columns indexed by $(1,2), (1,4), ...,
(1,g-1),(2,k+1)$ and by the last $g-2$ rows has rank $g-2$. Denote
by $(\lambda_1,...,\lambda_{g-2})$ the coordinates of the vector
given by the column of $Z_1$ indexed by $(1,3)$ and the last $g-2$
rows, with respect to the basis of ${\Co}^{g-2}$ given by the
columns of $C$. Then the coordinates of the columns of $A$ with
respect to the basis given by the columns of $C$ are given by the
following matrix which we will show to have maximal rank:
$${\small \left( \begin{array}{ccccccccccc}
\lambda_1-1&-\frac{1}{2}&-\frac{1}{3}&...&-\frac{1}{k-2}&0&...&0&\alpha_1&\beta_1\\
\lambda_2&\frac{1}{2}&0&...&0&0&...&0&0&0 \\
\lambda_3&0&\frac{1}{3}&...&0&0&...&0&0&0 \\
.&...&.&.&...&.&.&.&.&.& \\
.&...&.&.&...&.&.&.&.&.& \\
.&0&.&.&...&.&.&0&0&0& \\
\lambda_{k-2}&0&0&...&\frac{1}{k-2}&0&...&0&-\alpha_1&-\beta_1 \\
\lambda_{k-1}&0&0&...&0&\mu_{1,k+2}&...&\mu_{1,g-1}&\alpha_3&\beta_3 \\
\lambda_{k}&0&0&...&0&\mu_{2,k+2}&...&0&0&0 \\
&...&.&.&...&0&.&.&.&.& \\
.&...&.&.&...&.&.&0&.&0& \\
\lambda_{g-3}&0&0&0&...&0&0&\mu_{2,g-1}&0&\beta_4& \\
\lambda_{g-2}&0&0&0&...&\mu_{3,k+2}&..&\mu_{3,g-1}&\alpha_4&\beta_5& \\
\end{array} \right),}$$
where, for $j=k+2,...,g-1$,
$$\mu_{1,j} = \frac{1}{j-2} (\frac{kc_{j,1}}{c_{k+1,1}}-\frac{(k+1)c_{j,2}}{c_{k+1,2}}),
\ \mu_{2,j} = \frac{1}{j-2},$$
$$\mu_{3,j} = (\frac{k-1}{j-2}) (-\frac{c_{j,1}}{c_{k+1,1}}+\frac{c_{j,2}}{c_{k+1,2}}), \ \alpha_3= \frac{k-2}{2}, $$
$$\alpha_4= \frac{(k-1)(4-k^2)}{2k(k+1)}, \  \beta_3 = \frac{1}{g-k-1} (\frac{-3kc_{g-1,2}}{2c_{k+1,2}}+\frac{2(k-1)c_{g-1,1}}{c_{k+1,1}}),$$
$$\beta_4= \frac{1}{g-k-1}, \ \beta_5 =  \frac{1}{g-k-1} (\frac{3k(k-1)c_{g-1,2}}{2(k+1)c_{k+1,2}}-\frac{2(k-1)^2c_{g-1,1}}{kc_{k+1,1}}).$$

Substracting from each of the last two columns a suitable multiple of the $(k-2)$'s column, we can assume that $\alpha_1=\beta_1 =0$, hence the submatrix formed by the last $g-k$ columns and the first $k-2$ rows is zero.
The determinant of the submatrix given by the first $(k-2)$ rows and columns is $\frac{1}{(k-2)!} (\sum_{i=1}^{k-2} \lambda_i -1)$, and the determinant of the submatrix given by the last $g-k$ rows and columns is a non zero multiple of $$det\left( \begin{array}{ccc}
\mu_{1,g-1}&\alpha_3&\beta_3 \\
\mu_{2,g-1}&0&\beta_4 \\
\mu_{3,g-1}&\alpha_4&\beta_5 \\
\end{array} \right)= \frac{-(k-1)(k-2)^2(g-k-2)}{2k(k+1)^2(k+2)(g-3)(g-k-1)} \neq 0.$$
So it remains to show that $\sum_{i=1}^{k-2} \lambda_i \neq 1$. To do this, it suffices to show that  the matrix obtained by adding the row $(1,...,1)$ to the submatrix of $Z_1$ formed by the columns indexed by $(1,2),(1,3),...,(1,g-1),(2,k+1)$ and the last $g-2$ rows has rank $g-1$. This can be easily seen with a procedure similar to the one used in Proposition \ref{quad1}.

 \qed

\section{Surjectivity}
In this section we will prove by induction on the genus the surjectivity of $\mu_A$ for a general Prym-canonical binary curve $(C,A)$ of genus $\geq 20$.
\subsection{The 2nd Gaussian map}

Let us first of all analyze in detail the map $\mu_A$ of \eqref{gauss} when $C= C_1 \cup C_2$  is a Prym-canonical binary curve embedded in ${\proj}^{g-2}$ by $\omega_C \otimes A$, where $A \in Pic^0(C)$ is nontrivial of order 2.

Since ${\omega_C}_{|C_i} = K_{C_i}(D_i)$ where $D_i$ is the divisor of nodes in $C_i$, we have
$$H^0(S^2(\Omega^1_C) \otimes \omega_C^{\otimes 2}) \cong T \oplus (\oplus_{i=1,2} H^0(C_i, K^{\otimes 4}_{C_i}(2D_i))),$$
where $T$ is the torsion of $S^2(\Omega^1_C)$, which is supported at the nodes (see lemma 2 of \cite{ccm}). In fact we have an exact sequence
$$0 \rightarrow T \rightarrow S^2(\Omega^1_C) \rightarrow {\mathcal F}_C \rightarrow 0,$$
where $ {\mathcal F}_C$ is a non-locally free, rank 1, torsion free sheaf on $C$.

To prove the surjectivity of $\mu_A$ we will show the surjectivity of the components of $\mu_A$ on both non torsion and torsion parts of  $H^0(S^2(\Omega^1_C) \otimes \omega_C^{\otimes 2})$.

Consider first the non torsion component $\nu = \nu_1 \oplus \nu_2$, where
$$\nu_k: I_2(C)\rightarrow  H^0(C_k, K^{\otimes 4}_{C_k}(2D_k)) \cong H^0({\proj}^1, {\mathcal O}_{{\proj}^1}(2g-6)), \ k =1,2.$$

Recall that the curves $C_l$, $l=1,2$ are the images of the maps $\phi_l$ defined in \eqref{phi},
$\phi_l:{\proj}^1 \rightarrow {\proj}^{g-2}, \ l=1,2$
$$\phi_l(t,u):= (f_{1,l}(t),...,f_{g-1,l}(t)), \ f_{i,l}(t)= M_l(t)\frac{(\de_{i,l}t-c_{i,l})}{(t-a_{i,l})}.$$

Assume $Q \in I_2(C)$ is of the form (\ref{Qua}) where $s_{ij}$ are
solutions of (\ref{Qsistem}). Then using the local expression given in \eqref{mulocal}, we have
\begin{equation}\label{nuk}
\nu_k(Q) = \sum_{1\leq i<j\leq g-1} M^2_k(t)(\frac{\de_{i,k}t-c_{i,k}}{t-a_{i,k}})^{'}( \frac{\de_{j,k}t-c_{j,k}}{t-a_{j,k}})^{'}s_{ij}(dt)^4, \ k=1,2
\end{equation}

As an element of $H^0({\proj}^1, {\mathcal O}_{{\proj}^1}(2g-6))$, $\nu_k(Q)$ can be identified with the polynomial of degree $2g-6$ in $t$,
\begin{equation}\label{Rk} R_k(t) = \sum_{1\leq i<j\leq g-1} M^2_k(t)\frac{(c_{i,k}-\de_{i,k}a_{i,k})}{(t-a_{i,k})^2}
\frac{(c_{j,k}-\de_{j,k}a_{j,k})}{(t-a_{j,k})^2}s_{ij}.\end{equation}

To study the torsion component, we consider as in \cite{ccm} the restriction $\tau$ of $\mu_A$ to $ker(\nu)$, which lands in the torsion part $T$ of $H^0(S^2(\Omega^1_C) \otimes \omega_C^{\otimes 2})$.Then using Lemma 2 of \cite{ccm} one sees that the composition of $\tau$ with the projection on the torsion part $T_{P_h}$ at the nodes $P_1,...,P_{g-1},P_g$ is as follows:
a quadric $Q \in ker(\nu)$ as in (\ref{Qua}) is mapped to
\begin{equation}\label{torsion}
\begin{gathered}
dx dy \sum_{i \neq j}s_{ij}
f_{i,1}^{'}(a_{h,1})f_{j,2}^{'}(a_{h,2}) + 2xdx dy \sum_{i \neq
j}s_{ij}
f_{i,1}^{''}(a_{h,1})f_{j,2}^{'}(a_{h,2}) + \\
2y dx dy \sum_{i \neq j}s_{ij}
f_{i,1}^{'}(a_{h,1})f_{j,2}^{''}(a_{h,2}),
\end{gathered}
\end{equation}
where $h=1,...,g$, $s_{ij} = s_{ji}$ and $x$, $y$ are local coordinates around $P_h$ such that $C_1$ is given locally by $x=0$ and $C_2$ by $y=0$ and since $P_g$ is the image of $[0,1]$, we set $a_{g,1} = a_{g,2} =0$.
The description of the torsion at the point $P_{g+1}$ is similar:
\begin{equation}
\label{torsionPg+1}
\begin{gathered}
dx dy \sum_{i \neq j}s_{ij} (\de_{i,1}a_{i,1} - c_{i,1})(\de_{j,2}a_{j,2}-c_{j,2}) +\\ 2xdx dy \sum_{i \neq j}s_{ij} a_{i,1}(\de_{i,1}a_{i,1}-c_{i,1})(\de_{j,2}a_{j,2} - c_{j,2})  +\\
2y dx dy \sum_{i \neq j}s_{ij} a_{j,2}(\de_{i,1}a_{i,1} - c_{i,1})(\de_{j,2}a_{j,2} - c_{j,2}),
\end{gathered}
\end{equation}
where $s_{ij} = s_{ji}$ and $x$, $y$ are local coordinates around $P_{g+1}$ such that $C_1$ is given locally by $x=0$ and $C_2$ by $y=0$.

\subsection{Proof of surjectivity}
Let $C \subset {\proj}^{g-2}$ be a Prym-canonical binary curve embedded by $\omega_C \otimes A$, with $A^{\otimes 2} \cong \OO_C$, as in (\ref{deltac}) and set $k :=[\frac{g}{2}]$. Denote by $\tilde{C}$ the partial normalization of $C$ at the node
$P$, where  $P=P_k$ if $g =2k$, $P= P_{k+1}$ if $g = 2k+1$ and
by $p_1,p_2$ the preimages of $P$ in $\tilde{C}$. Observe that for a general choice of the $a_{i,j}$'s, the
projection $\pi$ from $P$ sends the curve $C$ to the Prym-canonical
model of $\tilde{C}$ in ${\proj}^{g-3}$ given by the line bundle
$K_{\tilde{C}} \otimes  A'$ where $A'$
corresponds to the point $(h'_1,...,h'_{g-1}, 1) \in {{\Co}^*}^{g}/{\Co}^*$, with $h'_i = 1$
for $i\leq [\frac{g-1}{2}]$, $h'_i = -1$ for
$i=[\frac{g-1}{2}]+1,...,g-1$,  as in section 2. In
fact if $g =2k$, we have $k-1 = [\frac{g-1}{2}]=:k'$ and $(\tilde{C},A')$ is  as in \eqref{phi}, \eqref{deltac} with
$a'_{i,j} = a_{i,j}$ for $i \leq k'$, $j =1,2$,  $a'_{i,j} =
a_{i+1,j}$ for $i \geq k'+1$, $j =1,2$.  If $g = 2k+1$ we have
$[\frac{g-1}{2}] = k$, so $(\tilde{C},A')$ is parametrized by
$a'_{i,j} = a_{i,j}$ for $i
\leq k$, $j =1,2$,  $a'_{i,j} = a_{i+1,j}$ for $i \geq k+1$,
$j=1,2$.

Consider the following commutative diagrams with horizontal exact
sequences
\begin{equation}
\label{chi} \xymatrix{
&0\ar[r] & \oplus_{i=1,2} H^0(C_i, K^{\otimes 4}_{C_i}(2\tilde{D_i}))\ar[r] & \oplus_{i=1,2} H^0(C_i, K^{\otimes 4}_{C_i}(2D_i)) \ar[r]  & \oplus_{i=1,2}\OO_{2p_i}\\
&0 \ar[r] &I_2({\tilde C}) \ar[u]^{\tilde\nu} \ar[r]& I_2(C)
\ar[u]^{\nu}\ar[ur]_{\chi} & }\end{equation}

\begin{equation}
\label{tau} \xymatrix{
&0\ar[r] & \tilde T \ar[r] & T  \ar[r] & T_P \\
&0 \ar[r] &\ker(\tilde\nu) \ar[u]^{\tilde\tau}  \ar[r]& ker(\nu)
\ar[u]^{\tau} \ar[ur]_{\tau_P} & }\end{equation}

where $D_i$ is the divisor of nodes of $C$ on $C_i$ and $\tilde{D_i} = D_i + p_i$ and $\nu, \tau$ and $\tilde\nu, \tilde\tau$ are the maps defined in the previous section for $C$ and $\tilde C$.
Hence, if  $\tilde\nu$ and $\chi$
($\tilde\tau$ and $\tau_P$, resp.) are surjective, then $\nu$
($\tau$, resp.) is also surjective.

\begin{TEO}\label{thm20}
If $(C=C_1\cup C_2, \omega_C\otimes A)$ is a Prym-canonical general
binary curve of genus $g \geq 20$, then $\mu_A$ is surjective for
$C$.
\end{TEO}
\proof The case $g=20$ is done by a direct computation with
Maple (see Appendix A). We then proceed by induction on $g$: the
commutativity of the diagrams \eqref{chi} \eqref{tau}  shows that
it is enough to prove the surjectivity of $\chi$ and $\tau_P$,
where $P=P_k$ for $g=2k$ and  $P=P_{k+1}$ for $g=2k+1$, as above.
Recall that the map $\nu$ is $\nu_1\oplus\nu_2$ where $\nu_1$ and
$\nu_2$ are defined in (\ref{nuk}), so we can write
$\chi=\chi_1\oplus\chi_2$, where $\chi_l$ is the composition of
$\nu_l$ with the restriction to $\OO_{2p_l}$, $l=1,2$. We want to
compute $\chi(Q)$, where $Q\in I_2(C)$. Notice that if $Q\in
I_2(C)$ is of the form \eqref{Qua}, with the $s_{ij}$'s satisfying
\eqref{qqij}, then $\chi_l(Q)$ is the pair
$(R_l(a_{r,l}),R'_l(a_{r,l}))$ (where $r=k$ for $g$ even and
 $r=k+1$ for $g$ odd) corresponding to the evaluation of the polynomial $R_l(t)$ ($l=1,2$) of \eqref{Rk}
 and of its derivative  at $P$. Recall that $R_l(t)$ is linear in the $s_{ij}$'s and
denote by $R^l_{i,j}(t)$ the coefficient of $s_{ij}$ in $R_l(t)$.

To prove the surjectivity  of $\chi$ we have to show that the
matrix $Y$ of size $(2g+2)\times \binom{g-1}{2}$ obtained by
concatenating vertically the matrix $Z$ in the proof of
Proposition \ref{Z}, and the matrix of size $4\times
\binom{g-1}{2}$  whose rows are the evaluations in $P$ of
$R^1_{i,j},(R^1_{i,j})',R^2_{i,j},(R^2_{i,j})'$ is of maximal
rank.

By formula \eqref{Rk} we have:
\begin{equation}
R^l_{i,j}(t)=(c_{i,l}-a_{i,l}\delta_{i,l})(c_{j,l}-a_{j,l}\delta_{j,l})\cdot\prod_{r\neq
i,j}(t-a_{r,l})^2. \end{equation} Therefore, if $i,j\neq n$,
$R^l_{i,j}(a_{n,l})=0$ and $(R^l_{i,j})'(a_{n,l})=0$. So it
remains to compute  $R^l_{i,k}(a_{k,l}),R^l_{k,j}(a_{k,l}) $ and
$(R^l_{i,k})'(a_{k,l}),(R^l_{k,j})'(a_{k,l})$ , for $g=2k$, and
$R^l_{i,k+1}(a_{k+1,l}),$ $R^l_{k+1,j}(a_{k+1,l}) $ and
$(R^l_{i,k+1})'(a_{k+1,l}),$ $(R^l_{k+1,j})'(a_{k+1,l})$, for
$g=2k+1$. If $g=2k$ and we denote by $D_{k,l}:=\prod_{r\neq
k}(a_{k,l}-a_{r,l})^2$, we have
$$ R^l_{i,k}(a_{k,l})= \frac{D_{k,l}\cdot
a_{i,l}a_{k,l}}{(a_{k,l}-a_{i,l})^2} ,\ R^l_{k,j}(a_{k,l})=-
\frac{D_{k,l}\cdot c_{j,l}a_{k,l}}{(a_{k,l}-a_{j,l})^2}
$$

$$ (R^l_{i,k})'(a_{k,l})= \frac{2D_{k,l}\cdot
a_{i,l}a_{k,l}}{(a_{k,l}-a_{i,l})^2}\cdot \sum_{r\neq
i,k}\frac{1}{(a_{k,l}-a_{r,l})}, $$
$$
(R^l_{k,j})'(a_{k,l})=- \frac{2D_{k,l}\cdot
c_{j,l}a_{k,l}}{(a_{k,l}-a_{j,l})^2}\cdot \sum_{r\neq
j,k}\frac{1}{(a_{k,l}-a_{r,l})}
$$

If $g=2k+1$, and we denote by $D_{k+1,l}:=\prod_{r\neq
k+1}(a_{k+1,l}-a_{r,l})^2$ we have
$$ R^l_{i,k+1}(a_{k+1,l})=-
\frac{D_{k+1,l}\cdot a_{i,l}c_{k+1,l}}{(a_{k+1,l}-a_{i,l})^2},\
R^l_{k+1,j}(a_{k+1,l})= \frac{D_{k+1,l}\cdot
c_{j,l}c_{k+1,l}}{(a_{k+1,l}-a_{j,l})^2}
$$

$$ (R^l_{i,k+1})'(a_{k+1,l})=- \frac{2D_{k+1,l}\cdot
a_{i,l}c_{k+1,l}}{(a_{k+1,l}-a_{i,l})^2}\cdot \sum_{r\neq
i,k+1}\frac{1}{(a_{k+1,l}-a_{r,l})}
$$
$$(R^l_{k+1,j})'(a_{k+1,l})= \frac{2D_{k+1,l}\cdot
c_{j,l}c_{k+1,l}}{(a_{k+1,l}-a_{j,l})^2}\cdot \sum_{r\neq
j,k+1}\frac{1}{(a_{k+1,l}-a_{r,l})}
$$

To show that the matrix $Y$ has maximal rank  $2g+2$ we will show
that the minor $det N$ is non zero, where $N$ is determined by the
columns indexed by $(1,i)$, $(2,j)$,
 with $2\leq i\leq g-1$, $3\leq j\leq g-1$, $(k, k+1), (k, g-1),
(g-2,g-1)$ and we choose the columns
$(3,k),(4,k),(k+1,g-2),(k+1,g-1)$, in the case $g=2k$, and the
columns $(3,k),(4,k+1),(k+1,g-2),(k+1,g-1)$, in the case $g=2k+1$.
Notice that the square submatrix of $N$ given by the first $2g-2$
rows and columns is the submatrix $Z_1$ of $Z$ introduced in
Prop.\ref{Z}, which is shown to be non singular for a general
choice of the $a_{i,l}$'s. The columns of the submatrix $G$ of $N$
given by its last four columns and its first $2g-2$ rows are
clearly also columns of $Z$ hence linearly dependent on the
columns of $Z_1$. Therefore we perform operations on the last four
columns of $N$ to bring $G$ to the zero matrix. So it suffices to
prove that the submatrix  $A$  of order 4, given by the last 4
rows and columns is nonsingular for general $a_{i,l}$. To do this
we choose the set of the $a_{i,l}$'s as in \eqref{scelta}, we
compute with Maple the determinant of $A$ and we see that as a
function of $k$ it does not vanish for any integer $k \geq 10$ (see
Appendix B). This proves that $\chi$ is surjective.

It remains to show that $\tau_P$ is surjective. Recall that
$\ker(\nu)$ is defined in $ I_2(C)$ by the vanishing of the
polynomials $R_l(t)$, $l=1,2$. By the description of the torsion
at the point $P$  given in \eqref{torsion}, we need to show the
rank maximality of the matrix $X$ of size $(2g+5)\times
\binom{g-1}{2}$ obtained by concatenating vertically the above
matrix $Y$ and the matrix of size $3\times \binom{g-1}{2}$ whose
rows are, for $g=2k$ (hence $P=P_k$)
$$(T_1)_{ij}=f_{i,1}^{'}(a_{k,1})f_{j,2}^{'}(a_{k,2})+f_{j,1}^{'}(a_{k,1})f_{i,2}^{'}(a_{k,2}), $$
$$(T_2)_{ij}=f_{i,1}^{''}(a_{k,1})f_{j,2}^{'}(a_{k,2})+  f_{j,1}^{''}(a_{k,1})f_{i,2}^{'}(a_{k,2}),$$
$$(T_3)_{ij}=f_{i,1}^{'}(a_{k,1})f_{j,2}^{''}(a_{k,2}) +  f_{j,1}^{'}(a_{k,1})f_{i,2}^{''}(a_{k,2})$$ and for
$g=2k+1$, hence $P=P_{k+1}$,
$$(T_1)_{ij}=f_{i,1}^{'}(a_{k+1,1})f_{j,2}^{'}(a_{k+1,2})+ f_{j,1}^{'}(a_{k+1,1})f_{i,2}^{'}(a_{k+1,2}),$$
$$ (T_2)_{ij}=f_{i,1}^{''}(a_{k+1,1})f_{j,2}^{'}(a_{k+1,2}) + f_{j,1}^{''}(a_{k+1,1})f_{i,2}^{'}(a_{k+1,2}),$$
$$(T_3)_{ij}=f_{i,1}^{'}(a_{k+1,1})f_{j,2}^{''}(a_{k+1,2}) + f_{j,1}^{'}(a_{k+1,1})f_{i,2}^{''}(a_{k+1,2}).$$

 We claim that the minor
$det M$ of the submatrix $M$ of $X$ determined by the $2g+5$
columns, indexed as the columns of $N$ plus $(5,k+1), (k,g-4),(k+1,g-3)$ if $g = 2k$, and
$(5,k+1), (k+1,g-4),(k+1,g-3)$ if $g = 2k+1$ is nonzero.
 This will conclude the proof that $\tau_P$ is surjective, hence the proof of the theorem.

As above the square submatrix of $M$ given by the first $2g-2$
rows and columns is the submatrix $Z_1$ of $Z$ introduced in
Prop.\ref{Z}, which is non singular for a general choice of the
$a_{i,l}$'s. The columns of the submatrix $H$
 of $M$ given by its last seven  columns and its first $2g-2$
rows are clearly also columns of $Z$ hence linearly dependent on
the columns of $Z_1$. Therefore we perform operations on the last
seven columns of $M$ to bring $H$ to the zero matrix. So it
suffices to prove that the submatrix  $B$ of order 7, given by the
last 7 rows and columns is nonsingular for general $a_{i,l}$. To
this purpose we choose the set of the $a_{i,l}$'s as in
\eqref{scelta}, we compute again with Maple the determinant of $B$
and we see that for any integer $k \geq 10$ it does not vanish (see Appendix B).
This proves that $\tau_P$ is surjective, hence by induction
$\mu_A$ is surjective. \qed

\section{The class}

In the previous section we have proved by semicontinuity that the
2nd Gaussian map $\mu_A: I_2(C) \rightarrow H^0(S^2(\Omega^1_C)
\otimes K_C^{\otimes 2})$ has maximal rank for the general pair
$[C,A]$ in ${\mathcal R}_{20}$. Notice that for $g =20$,
$dim(I_2(C)) = dim(H^0(S^2(\Omega^1_C) \otimes K_C^{\otimes 2})) =
133. $ Consider the locus ${\mathcal D} = \{[C,A] \in {\mathcal
R}_{20} \ | \ rk(\mu_A) <133 \}$. We have proved that ${\mathcal
D} \neq {\mathcal R}_{20}$, hence, if it is not empty, it is an
effective divisor in ${\mathcal R}_{20}.$ Let $\pi:
\overline{\mathcal R}_{g} \rightarrow \overline{{\mathcal M}}_g$
be the finite map which extends the forgetful map ${\mathcal R}_g
\rightarrow {\mathcal M}_g$ (see \cite{fl} Section 1). The partial
compactification $\tilde{\mathcal R}_g$ of ${\mathcal R}_g$
introduced in \cite{fl} Section 1 is the inverse image
$\pi^{-1}(\tilde{\mathcal M}_g)$, where $\tilde{\mathcal
M}_g:={\mathcal M}_g\cup \tilde{\Delta }_0$ and $\tilde{\Delta
}_0$ is the locus of one-nodal irreducible curves. Denote by $f:
{\mathcal X} \rightarrow \tilde{{\bf{R}}}_g$ the universal family
and by ${\mathcal P} \in Pic({\mathcal X})$ the corresponding Prym
bundle as in \cite{fl} 1.1. Assume $g =20$, then if
$\tilde{{\mathcal D}}$ is the closure of ${\mathcal D}$ in
${\tilde{\mathcal R}_{20}}$, $\tilde{{\mathcal D}}$ is the
degeneracy locus of the map
$$\tilde{\mu}: \tilde{{\mathcal I}}_2 \rightarrow f_*( (\omega_f\otimes {\mathcal P})^{\otimes 2}\otimes
{S^2(\Omega^1_f)})\cong f_*(\omega_f^{\otimes 4}\otimes {\mathcal
P}^{\otimes 2}\otimes {\mathcal I}_Z^{\otimes 2}),$$ of
\eqref{mutilde}. Denote by ${\mathcal F}_i :=
f_{*}(\omega_f^{\otimes i} \otimes {\mathcal P}^{\otimes i})$.
Using Grothendieck-Riemann-Roch and Proposition 1.6 of \cite{fl}
one computes as in Proposition 1.7 of \cite{fl}
$$
 c_1({\mathcal F}_i) = \frac{i(i-1)}{2} (12 \lambda - \delta'_0- \delta''_0 - 2 \delta_0^{ram}) +
 \lambda - \frac{i^2}{4} \delta_0^{ram},
$$
where $\lambda$ is the pullback of the Hodge class $\lambda \in {\overline{\mathcal M}}_g$ and $ \delta'_0$, $\delta''_0$, and $\delta_0^{ram}$ are the boundary classes defined in \cite{fl} section 1.
So we have
$$
 c_1({\mathcal F}_1) = \lambda - \frac{\delta_0^{ram}}{4}, \ \  c_1({\mathcal F}_2) = 13 \lambda -
 \delta'_0- \delta''_0- 3\delta_0^{ram}, \ \
 c_1(S^2({\mathcal F}_1)) = 20 \cdot c_1({\mathcal F}_1)= 20 \lambda - 5 \delta_0^{ram},
$$
therefore
$$c_1({\mathcal I}_2) = c_1(S^2{\mathcal F}_1) - c_1({\mathcal F}_2) = 7 \lambda + \delta'_0 +
\delta''_0 - 2 \delta_0^{ram}.$$

Notice that  by Grothendieck-Riemann-Roch we have
 $$c_1( f_*( \omega_f^{\otimes 4} \otimes {\mathcal P}^{\otimes 2} \otimes {\mathcal I}_Z^{\otimes 2})) =
 f_*[(1 + c_1( \omega_f^{\otimes 4} \otimes {\mathcal P}^{\otimes 2})
 + \frac{1}{2} c_1^2( \omega_f^{\otimes 4} \otimes {\mathcal P}^{\otimes 2})-2 [Z])
 \cdot (1 - \frac{c_1(\omega_f)}{2} + \frac{ c_1(\omega_f)^2+ [Z]}{12}) ]_2$$
$$= 73 \lambda -8 (\delta'_0 + \delta''_0) -17 \delta_0^{ram}, $$
since $f_*(c_1(\omega_f) \cdot {\mathcal P}) =0$, $f_*(c_1({\mathcal P})^2) = - \delta_0^{ram}/2$,
by Proposition 1.6 of \cite{fl} and by Mumford's formula, $f_*(c_1(\omega_f)^2) = 12 \lambda + f_*([Z])$
and $f_*[Z] = \delta'_0 + \delta''_0+ 2 \delta_0^{ram}$ (\cite{fl}, 1.1). So, finally we have
$$ c_1(\tilde{{\mathcal D}}) = c_1( f_*( \omega_f^{\otimes 4} \otimes {\mathcal P}^{\otimes 2}
\otimes {\mathcal I}_Z^{\otimes 2}))\cdot rk({\mathcal I}_2 )- c_1({\mathcal I}_2)\cdot
rk  (f_*( \omega_f^{\otimes 4} \otimes {\mathcal P}^{\otimes 2}\otimes {\mathcal I}_Z^{\otimes 2}))= $$
$$=133 (66 \lambda -9 (\delta'_0 + \delta''_0) -15 \delta_0^{ram}),$$
and $c_1({\mathcal D}) = 8778 \lambda$, hence ${\mathcal D}$ is an effective divisor in
${\mathcal R}_{20}$, $\tilde{{\mathcal D}}$ is an effective divisor in $\tilde{{\mathcal R}}_{20}$
and if we denote by $\overline{{\mathcal D}}$ the closure of ${\mathcal D}$ in $\overline{\mathcal R}_{20}$,
we have computed
\begin{equation}
\label{Dbar}
c_1(\overline{{\mathcal D}}) = 133 (66 \lambda -9 (\delta'_0 + \delta''_0) -15 \delta_0^{ram}- ...)
\end{equation}
In fact, since the partial compactification $\tilde{{\mathcal
R}}_{g} \subset \overline{\mathcal R}_{g}$ has the property that
$\pi^{-1}({\mathcal M}_g \cup \Delta_0)-\tilde{{\mathcal R}}_{g} $
has codimension $\geq 2$, the expression (\ref{Dbar}) computes the
coefficients of $\lambda$, $\delta'_0$, $\delta''_0$,
$\delta_0^{ram}$ in $c_1(\overline{{\mathcal D}})$.

\begin{REM}

\begin{itemize}
\item Using proposition 1.9 of \cite{fl} one can find lower bounds on some of the other boundary coefficients of $\overline{{\mathcal D}}$.
\item Pushing forward $\overline{{\mathcal D}}$, one gets
$$c_1(\pi_*(\overline{{\mathcal D}})) =
133 (66(2^{40}-1) \lambda -(33\cdot2^{38}-9) \delta_0-...),$$
hence its slope is $\geq 8+\frac{2}{3023656976381}$.
\end{itemize}
\end{REM}

\appendix

\section*{Appendix: Maple scripts for computations}   

{\bf{A. Surjectivity for $g=20$}} We list here the Maple script we
run. We will explain it afterwards: for this purpose, we added
line numbers.

\bigskip\noindent{\small
\verb!a[1]:=[25,35,54,47,67,97,73,81,22,33,76,27,38,44,58,69,63,80,99]:!
\\\verb!a[2]:=[1,3,4,5,6,7,8,9,10,11,12,13,14,15,16,17,18,19,20]:!
\\\verb!listsij:=[seq(seq(s[i,j],j=i+1..19),i=1..19)]:!
\\\verb!t[1]:=[0,0,0,0,0,0,0,0,0,0,1,1,1,1,1,1,1,1,1]:!
\\\texttt{\hspace*{-4ex}5\settowidth{\spazio}{5}\hspace{-\spazio}\hspace{4ex}}%
\verb!A1:=mul(a[1][i],i=1..19): !
\\\verb!A2:=mul(a[2][i],i=1..19): !
\\\verb!d2:=1:!
\\\verb!d1:= -d2*A1/A2:!
\\\verb!delta:=[1,1,1,1,1,1,1,1,1,1,0,0,0,0,0,0,0,0,0]:!
\\\texttt{\hspace*{-5ex}10\settowidth{\spazio}{15}\hspace{-\spazio}\hspace{5ex}}%
\verb!c[1]:=[seq(t[1][i]*a[1][i]*d1/A1,i=1..19)]:!
\\\verb!c[2]:=[seq(t[1][i]*a[2][i]*d2/A2,i=1..19)]:!
\\\verb!Z:=Matrix([seq([seq(seq((delta[i]*delta[j])*add(a[1][i]^m*a[1][j]^(h-m),m=0..h)-!
\\\verb!(delta[j]*c[1][i]+delta[i]*c[1][j])*add(a[1][i]^m*a[1][j]^(h-1-m),m=0..h-1)+ !
\\\verb!(c[1][i]*c[1][j])*add(a[1][i]^m*a[1][j]^(h-2-m),m=0..h-2), j=i+1..19),i=1..19)],!
\\\verb! h=0..19), seq([seq(seq((delta[i]*delta[j])*add(a[2][i]^m*a[2][j]^(h-m),m=0..h)-!
\\\verb!(delta[j]*c[2][i]+delta[i]*c[2][j])*add(a[2][i]^m*a[2][j]^(h-1-m),m=0..h-1)+!
\\\verb! (c[2][i]*c[2][j])*add(a[2][i]^m*a[2][j]^(h-2-m),m=0..h-2), j=i+1..19),i=1..19)],!
\\\verb!h=0..19)]):!
\\\texttt{\hspace*{-5ex}19\settowidth{\spazio}{15}\hspace{-\spazio}\hspace{5ex}}%
\verb!Zref:=Gausselim(Z,'r0') mod 131:!
\\\verb! r0;!
           \\\verb!                     38!
\\\verb!M[1] := mul(t-a[1][i], i = 1 .. 19):!
\\\verb!M[2] := mul(t-a[2][i], i = 1 .. 19):!
\\\verb!for i from 1 to 19 do phi1[1, i] := diff(M[1]*(delta[i]*t-c[1][i])/(t-a[1][i]), t):!
\\\texttt{\hspace*{-5ex}25\settowidth{\spazio}{5}\hspace{-\spazio}\hspace{5ex}}%
\verb!phi1[2, i] := diff(M[2]*(delta[i]*t-c[2][i])/(t-a[2][i]), t) end do:!
\\\verb!R[1] := add(add(s[i, j]*phi1[1, i]*phi1[1, j], j = i+1 .. 19), i = 1 .. 19):!
\\\verb!R[2] := add(add(s[i, j]*phi1[2, i]*phi1[2, j], j = i+1 .. 19), i = 1 .. 19):!
\\\verb!Eqskernu := [seq(seq(primpart(coeff(R[k], t, n)), n = 0 .. 34), k = 1 .. 2)]:!
 \\\verb!K:= Gausselim(linalg[stackmatrix](Zref,linalg[genmatrix](Eqskernu,listsij)),'r1')!
\\\verb!mod 131):!
\\\verb! r1;!
  \\\verb!              108!
\\\texttt{\hspace*{-5ex}32\settowidth{\spazio}{15}\hspace{-\spazio}\hspace{5ex}}%
\verb!for i from 1 to 19 do phi2[1,i]:= diff(phi1[1,i],t): phi2[2,i] := diff(phi1[2,i],t): !
\\\verb!phi1e0[1,i]:= eval(phi1[1,i], t = 0): phi2e0[1,i] := eval(phi2[1,i], t = 0): !
\\\verb!phi1e0[2,i] := eval(phi1[2,i], t = 0): phi2e0[2,i]:= eval(phi2[2,i], t = 0): !
\\\verb!for h to 19 do phi1e[1, i, h]:= eval(phi1[1,i], t = a[1][h]):!
\\\verb!phi2e[1, i, h]:= eval(phi2[1,i], t = a[1][h]):!
\\\verb!phi1e[2, i, h]:= eval(phi1[2,i], t = a[2][h]):!
\\\verb!phi2e[2, i, h] := eval(phi2[2,i], t = a[2][h]) end do end do:!
\\\texttt{\hspace*{-5ex}39\settowidth{\spazio}{15}\hspace{-\spazio}\hspace{5ex}}%
\verb!for h from 1 to 19 do!
\\\verb!tors[h,1]:= add(add(s[i,j]*(phi1e[1,i,h]*phi1e[2, j, h]+!
\\\verb!phi1e[1,j,h]*phi1e[2,i,h]), j = i+1 .. 19), i = 1 .. 19):!
\\\verb!tors[h,2]:= add(add(s[i,j]*(phi2e[1,i,h]*phi1e[2, j, h]+!
\\\verb!phi2e[1,j,h]*phi1e[2,i,h]), j = i+1 .. 19), i = 1 .. 19):!
\\\verb!tors[h,3]:= add(add(s[i,j]*(phi1e[1,i,h]*phi2e[2, j, h]+!
\\\texttt{\hspace*{-5ex}45\settowidth{\spazio}{15}\hspace{-\spazio}\hspace{5ex}}%
\verb!phi1e[1,j,h]*phi2e[2,i,h]), j = i+1 .. 19), i = 1 .. 19) end do:!
\\\verb!tors[20,1]:= add(add(s[i j]*(phi1e0[1,i]*phi1e0[2, j]+phi1e0[1, j]*phi1e0[2, i]), !
\\\verb!j = i+1 .. 19), i = 1 .. 19): !
\\\verb!tors[20,2]:= add(add(s[i,j]*(phi2e0[1,i]*phi1e0[2, j]+phi2e0[1, j]*phi1e0[2, i]), !
\\\verb!j = i+1 .. 19), i = 1 .. 19):!
\\\verb!tors[20,3]:= add(add(s[i,j]*(phi1e0[1,i]*phi2e0[2, j]+phi1e0[1, j]*phi2e0[2, i]),!
\\\verb!j = i+1 .. 19), i = 1 .. 19):!
\\\texttt{\hspace*{-5ex}52\settowidth{\spazio}{15}\hspace{-\spazio}\hspace{5ex}}%
\verb!tors[21,1]:= add(add(s[i,j]*((delta[i]*a[1][i]-c[1][i])*!
\\\verb!(delta[j]*a[2][j]-c[2][j])+!
\\\verb!(delta[j]*a[1][j]-c[1][j])*(delta[i]*a[2][i]-c[2][i])), j = 1 .. 19), i = 1 .. 19):!
\\\verb!tors[21,2]:= add(add(s[i, j]*((delta[i]*a[1][i]-c[1][i])*a[1][i]*!
\\\verb!(delta[j]*a[2][j]-c[2][j])+!
\\\verb!(delta[j]*a[1][j]-c[1][j])*(delta[i]*a[2][i]-c[2][i])*!
\\\verb!a[1][j]), j = 1 .. 19), i = 1 .. 19):!
\\\verb!tors[21,3]:= add(add(s[i, j]*((delta[i]*a[1][i]-c[1][i])*a[2][j]*!
\\\verb!(delta[j]*a[2][j]-c[2][j])+!
\\\verb!(delta[j]*a[1][j]-c[1][j])*(delta[i]*a[2][i]-c[2][i])*!
\\\verb!a[2][i]), j = 1 .. 19), i = 1 .. 19):!
\\\texttt{\hspace*{-5ex}63\settowidth{\spazio}{15}\hspace{-\spazio}\hspace{5ex}}%
\verb!Eqskertau:= [seq(seq(primpart(tors[h,l]), l = 1 .. 3), h = 1 .. 21)]:!
\\\verb!Gausselim(linalg[stackmatrix](K, linalg[genmatrix](Eqskertau, listsij)),'r2')!
\\\verb!mod 131):!
\\\verb!r2;!
\\\verb!        171!
\\\verb!r2-r1;!
\\\verb!                                63!
}

In lines 1--2, we define the $a_{i,j}$'s which will be used. We
chose them randomly.  In line 3, we collect the unknowns
$\{s_{i,j}\}_{1\le i < j \le g}$ in the list \texttt{listsij}:
there are $\binom{g-1}{2}$ of them.
In line 4,5,6 we define $\texttt{t[1]}$ which is the vector $P_g$ as in (\ref{lemmadeltac}), $A_i = \prod_{i=1,...,19} a_{r,i}$, $i=1,2$. In line 7,8 we define $d_2=1$, $d_1 = -\frac{d_2 A_1}{A_2}$, as in (\ref{deltac}). In line 9 we define the vector $\delta$ whose components are $\delta_{i,1} = \delta_{i,2}$, as in (\ref{deltac}).
In lines 10,11 we collect the $c_{i,1}$, $c_{i,2}$ $i=1,...,19$ as in (\ref{deltac}), and we call them $ \texttt{c[1][i]}$, $ \texttt{c[2][i]}$. These data give the curve $C$ and the line bundle $A$ as in (\ref{lemmadeltac}).
In lines 13--19, we define the matrix \texttt{Z} associated to the
linear system \eqref{qqij}, whose solutions give us the quadrics in
$I_2(C)$, cf. Proposition \ref{Z}. In line 20,
Maple computes the rank \texttt{r0} of \texttt{Z} via Gaussian
elimination, by calculating modulo 131 to speed up computations.
The resulting matrix  is called \texttt{Zref}.
As expected by Proposition \ref{Z},  Maple finds
$\texttt{r0}=38=2g-2$ and it prints it in line 21. In lines 22, 23 we define  \texttt{M[j]}  as $\prod_{i=1,...,19}(t-a_{i,j})$, $j=1,2$. In lines 24, 25 we define  \texttt{phi1[j, i]}  as the i-th component of $\frac{d}{dt}(\phi_j(t,1))$, $j=1,2$, where $\phi_j$ is defined in \eqref{phi}. In lines 26, 27 we define \texttt{R[k]} as the polynomial $R_k(t)$ of  \eqref{Rk}, $k=1,2$.
In line 28, we collect in \texttt{EqsKerNu} the list of equations
which determine $\ker(\nu)$,
cf.\ the definition of $\nu$ in \eqref{nuk}.
In lines 29--31, Maple computes the rank \texttt{r1} of the linear system
$\texttt{EqsKerNu} \cap\ker(\texttt{Zref})$,
again via Gaussian elimination modulo 131,
and the resulting matrix is called \texttt{K}.
Maple finds that $\texttt{r1}=108$ and it prints it in line 31. Since $rank(\nu) = \texttt{r1}- \texttt{r0} = 70 = 2 h^0(\OO_{\proj^1}(2g-6))$, for $g=20$, we have shown that $\nu$ has maximal rank.
In line 32, we define the 2nd derivative \texttt{phi2} of the $\phi_j(t,1)$'s of \eqref{phi}. In lines 33-34 we define the evaluations \texttt{phi1e0}, \texttt{phi2e0}, of the first and the second derivatives of  the $\phi_j(t,1)$'s at $t=0$, i.e. at the point $P_{20}$ and in lines 35-38 we define their evaluations at the points $P_i$, $i =1,...,19$.
Using them, in lines 39--51 we compute the torsion at $P_i$, $i=1,\ldots,20$,
cf.\ \eqref{torsion}, and, in lines 52--62,
the torsion at the point $P_{21}$, cf.\ \eqref{torsionPg+1}.
In lines 63 we collect in \texttt{EqsKerTau} the equations which
determine $\ker(\tau)$ and Maple computes the rank \texttt{r2}
of $\texttt{EqsKerTau} \cap\ker(\text{\texttt{K}})$,
via Gaussian elimination modulo 131 as before.
Maple finds that $\texttt{r2}=171$, therefore the rank of $\tau$ is
$\texttt{r2}-\texttt{r1}=171-108=63= 3(g+1)= dim(T)$, $(g=20)$, hence also $\tau$ has maximal rank.
So we have shown that $\mu_A$ is surjective for $g=20$.

{\bf{B. Results of computations in Thm.\ref{thm20}}}

Here we give the formulas of the determinants of the matrices $A$
and $B$ in the proof of Thm.\ref{thm20}. The Maple files of these
computations are available under request to the authors.

{\small If $g=2k$, $det A=\frac{ -4(4k^5+14k^4+15k^3+k^2-7k+1)}{p_1}$ where $$p_1=
   k^2 (k - 4)^3  (k + 1) (2 k - 1) (k - 1)^3 (k^2  - 4) (2 k^3  - 9 k^2
     + 12 k - 4) (k - 3)^2 ((2 k-1)! )^4 .$$

$det(B)=\frac{147456}{ 5  }(k-5)(k^4-9k^3+16k^2+3k-8)\cdot p_2\cdot p_3/(q_1 \cdot q_2),$ where
$$p_2=16k^9-14k^8-87k^7+121k^6+75k^5-138k^4-52k^3+54k^2+18k-12,$$
$$p_3=6k^5-k^4-12k^3-4k^2+12k-4$$
$$q_1 = k^3(k-3)^3(k^2-4k+4)(-3k+k^2+2)(2k-1)(k-1)^8a(k-4)^7(2k^3-9k^2+12k-4)$$
$$q_2= (k^2-4)(2k^2-7k+6)((2k)!)^8$$

and one can check all the functions appearing in these expressions
do not have any integral zero $k \geq 10$.

If $g=2k+1$, $ detA= \frac{-16}{15\cdot p_4}\cdot (k-4)p_5$ where
$$p_4=(-2+k)^4(2k-1)k^2(k-3)^3(k-1)^6(k+2)(k+1)((2+2k)!^2)),$$
$$p_5=(1168k^{14}+2216k^{13}-22360k^{12}-41218k^{11}+17145k^{10}+47730k^9+46525k^8+38736k^7$$
$$-70488k^6-58080k^5+35288k^4+14726k^3-6093k^2+66k-465).$$

$ det(B)=-\frac{3072}{25}(k-5) \cdot p_6 \cdot p_7/p_8,$ where
$$p_6=(270336+1257472k+25884500k^4-5217504k^2-15573704k^3-6492143k^{17}+68438542k^5$$
$$-28031103k^6-108784825k^7-49730235k^8-30298961k^9 +50987804k^{10}+197670424k^{11}+60883960k^{12}$$
$$-162484142k^{13}-9462204k^{18} -79976k^{19}+945456k^{20}+45632k^{21}-44288k^{22}-1216k^{23}+768k^{24}$$
$$-92612465k^{14}+54292657k^{15}+44402735k^{16}),$$
$$p_7=2+4k-k^2+(2k+1)!(k+1)(k^4-2k^3-k^2+12k+4),$$
$$p_8=ak^6(k-4)^5(k-1)(k^3-4k^2+5k-2)^3(2k^2-7k+3)(k^2-3k+2)^3(k+1)^3(k-3)^5(k+2)(2k-3)(4k^2-1)((2k+1)!)^7$$
and again one can check all these functions do not have any
integral zero $k \geq 10$.}

\end{document}